\def\version{21.2.2023}
\def\users{us}  %
\def\users{final-layout}   
\documentclass[12pt]{article}
\usepackage{upgreek}

\usepackage{xcolor}
\usepackage{bm,amsmath,amsthm,hyperref,amsfonts,amssymb,color}
\usepackage{mathrsfs} 
\usepackage{cite}

\usepackage{ifthen}
\ifthenelse{\equal{\users}{final-layout}}{}{
\usepackage{fancyhdr}
\pagestyle{fancy}
\headheight=28pt\headwidth=17cm
\definecolor{gray}{gray}{0.5}
\rhead{\color{gray}Fluid-solid interaction Eulerian\\
T.Roub\'\i\v cek}
\chead{}
\lhead{Version\,\version, file:\,\jobname.tex
\\
compiled:
\number\day.\number\month.\number\year\ at
\the\hour:\ifnum\minute<10 0\fi\the\minute\ h\ \ \ \ \ 
}
}

\definecolor{labelkey}{rgb}{1.,.2,0.}

\usepackage[normalem]{ulem}
\usepackage{ifthen}
\usepackage{color}

\ifthenelse{\equal{\users}{final-layout}}{

	\newcommand{\COMMENT}[1]{}
	\newcommand{\COMMENTGT}[1]{}
	\newcommand{\TODO}[1]{}
	\newcommand{\INTERNAL}[1]{}
	\newcommand{\QUESTION}[1]{}
	\newcommand{\DELETE}[1]{}

	\newcommand{\REM}[1]{\marginpar{\bfseries\tiny{\color{blue}}}}
    \newcommand{\MARGINOTE}[1]{}
}
{
	
	\newcommand{\COMMENT}[1]{{\color{red}\uuline{#1}\color{black}}}
	\newcommand{\COMMENTGT}[1]{{\hfill\large\color{red}***{#1}***\color{black}\hfill}\\}
	\newcommand{\TODO}[1]{{\color{red}\uuline{#1}\color{black}}}
	\newcommand{\INTERNAL}[1]{\footnote{#1}}
	\newcommand{\QUESTION}[1]{{\color{brown}\uuline{#1}\color{black}}}
	\newcommand{\DELETE}[1]{{\color{red}\sout{#1}\color{black}}}

	\newcommand{\REM}[1]{\marginpar{\bfseries\tiny{\color{blue}#1}}}
\newcommand{\MARGINOTE}[1]{\marginpar{\color{red}\tiny\texttt{#1}}}
}

\newcount\hour \newcount\minute
\hour=\time
\divide \hour by 60
\minute=\time
\loop \ifnum \minute > 59 \advance \minute by -60 \repeat

\newcommand\DT[1]{\mathchoice
                 {{\buildrel{\hspace*{.1em}\text{\LARGE.}}\over{#1}}}
                 {{\buildrel{\hspace*{.1em}\text{\LARGE.}}\over{#1}}}
                 {{\buildrel{\hspace*{.1em}\text{\Large.}}\over{#1}}}
                 {{\buildrel{\hspace*{.1em}\text{\large.}}\over{#1}}}}
\newcommand\pdt[1]{\frac{\partial{#1}}{\partial t}} 
\newcommand{\lineunder}[2]{\LU{\begin{array}[t]{c}\underbrace{#1}\vspace*{.5em}\end{array}}{\mbox{\footnotesize\rm #2}}}
\newcommand{\LU}[2]{\begin{array}[t]{c}#1\vspace*{-1em}\\_{#2}\end{array}}
\newcommand{\linesunder}[3]{\LSU{\begin{array}[t]{c}\underbrace{#1}\vspace*{.5em}\end{array}}{\mbox{\footnotesize\rm #2}}{\mbox{\footnotesize\rm #3}}}
\newcommand{\LSU}[3]{\begin{array}[t]{c}#1\vspace*{-1em}\\_{#2}\vspace*{-.5em}\\_{#3}\end{array}}

\newcommand{\Item}[2]{\parbox[t]{.055\textwidth}{#1}\hfill
\parbox[t]{.945\textwidth}{#2}\vspace*{.8mm}}

\newcommand{\divG}{\mathrm{div}_{\scriptscriptstyle{\hspace*{-.1em}\varGamma}}^{}}

\newcommand{\nablaG}{\nabla_{\scriptscriptstyle{\hspace*{-.1em}\varGamma}}^{}}
\def\Vdots{\!\mbox{\setlength{\unitlength}{1em}
\begin{picture}(0,0)
\put(-.07,0){.}
\put(-.07,.3){.}
\put(-.07,.6){.}
\end{picture}\hspace*{.2em}}}
\usepackage{mathrsfs}   
\usepackage{eucal}      

  \def\bbI{{\mathbb I}}

\def\FG{\boldsymbol}
   
 \def\ee{{\FG e}} \def\ff{{\FG f}}

 \def\nn{{\FG n}}  
  \def\rr{{\FG r}} 
 \def\tt{{\FG t}}  
\def\vv{{\FG v}}  \def\xx{{\FG x}} 
\def\yy{{\FG y}}  
\def\AA{{\FG A}} \def\BB{{\FG B}}  
\def\DD{{\FG D}} \def\EE{{\FG E}}
\def\FF{{\FG F}}

 \def\TT{{\FG T}} 
  \def\XX{{\FG X}}

\newcommand{\R}{\mathbb R}
\newcommand{\N}{\mathbb N}
\newcommand{\Nabla}{{\nabla}}
\newcommand{\pl}{\partial}
\newcommand{\eq}[1]{(\ref{#1})}
\renewcommand{\d}{\mathrm d}  
\newcommand{\barOmega}{\hspace*{.2em}{\overline{\hspace*{-.2em}\varOmega}}}

\newcommand{\jump}[1]{[\hspace*{-.15em}[#1]\hspace*{-.15em}]_{\varGamma_\text{\sc fs}(t)}^{}}
\newcommand{\Jump}[1]{\big[\hspace*{-.25em}\big[#1\big]\hspace*{-.25em}\big]_{\varGamma_\text{\sc fs}(t)}^{}}

\newcommand{\OmegaSt}{\varOmega_\text{\sc s}(t)}
\newcommand{\OmegaFt}{\varOmega_\text{\sc f}(t)}
\newcommand{\OmegaS}{\varOmega_\text{\sc s}}
\newcommand{\OmegaF}{\varOmega_\text{\sc f}}

\newcommand{\wt}{\widetilde}

\newtheorem{theorem}{Theorem}[section]
\newtheorem{lemma}[theorem]{Lemma}
\newtheorem{definition}[theorem]{Definition}

\newtheorem{proposition}[theorem]{Proposition}

\newtheorem{remark}[theorem]{Remark}

\newcommand{\TTRxi}{\TT_\text{\sc r}^{\bm\xi}}
\newcommand{\TTRexi}{\TT_{\text{\sc r},\LAM}^{\bm\xi}}
\newcommand{\TTRexik}{\TT_{\text{\sc r},\LAM}^{\bm\xi_k}}
\newcommand{\phiRxi}{\varphi_\text{\sc r}^{\bm\xi}}
\newcommand{\DISRxi}{\DIS_\text{\sc r}^{\bm\xi}}

\numberwithin{equation}{section}

\usepackage{graphicx}
\usepackage{psfrag} 
\newcounter{myfigure}
\newenvironment{my-picture}[3]{\refstepcounter{myfigure}\label{#3}\setlength{\unitlength}{1cm}\begin{picture}(#1,#2)}{\end{picture}}

\begin{document}

\allowdisplaybreaks

\pretolerance=10000  

\vspace*{4em}

\noindent
{\Large\bf Interaction of  finitely-strained viscoelastic multipolar\\[.3em]
solids and fluids by an Eulerian approach}

\bigskip\bigskip\bigskip

\noindent
{\large Tom\'a\v s Roub\'\i\v cek}

\bigskip\bigskip\bigskip

{\small
\noindent
{\bf Abstract}.
A mechanical interaction of compressible
viscoelastic fluids with viscoelastic solids in Kelvin-Voigt
rheology using the concept of higher-order  (so-called 2nd-grade
multipolar) viscosity is investigated in a quasistatic variant. The
no-slip contact between fluid and solid is considered and the
Eulerian-frame return-mapping technique is used for both the fluid and
the solid models, which allows for a ``monolithic'' formulation of this
fluid-structure interaction problem. Existence and a certain regularity
of weak solutions is proved by a Schauder fixed-point argument combined
with a suitable regularization.

\medskip


\medskip

\medskip

}

\bigskip

\baselineskip=16pt

\def\TRACTION{\bm{f}}
\def\GRAVITY{\bm{g}}
\def\rhoR{\varrho_\text{\sc r}^{}}
\def\rhoF{\varrho_\text{\sc f}^{}}
\def\rhoS{\varrho_\text{\sc s}^{}}
\def\rhoRxi{\varrho_\text{\sc r}^{\bm\xi}}
\def\rhoRxik{\varrho_\text{\sc r}^{\bm\xi_k}}
\def\LAM{\varepsilon}
\def\DIS{\DD}

\section{Introduction}

The mechanical interaction between fluids and solids, briefly a
{\it fluids-solid} (or {\it fluid-structure}) {\it interaction}
(FSI), is a largely scrutinized task in continuum mechanics. Numerically, there
are more than dozen of methods to cope with this problem, cf.\
\cite[Fig.1]{JaKaMa19CNEF}  for a survey. Here we focus on the so-called 
{\it Eulerian-frame Reference/Return map Technique} (RMT).
Actually, it can be seen as {philosophically}  similar to 
the popular arbitrary Lagrangian-Eulerian (ALE) formulation which
maps Eulerian fluid onto a fixed reference framework to be coupled with the
Lagrangian solid domain. Yet, RMT goes just opposite when mapping the solid
into actual Eulerian framework to be compatible with fluid, cf.\
\cite[Ch.6]{Rich17FSIM} or also \cite{WJWP20ESOM} for a comparison.

Engineering literature exploiting the return mapping (sometimes
confining on incompressible fluids or  on  a ``mushy'' interface by a
level-set method) as \cite{BeFoIo22EFVM,CPST17NSEM,CoMaMi08EFLS,Dunn06EAFS,DunRan06AFEA,DuRaRi10NSFS,Frei16EFEM,ISTT11IFEM,JaKaMa19CNEF,KaRyNa60RMTF,LaRBQu13FEFE,RanRic10AFEM,RichFEFF,RicWic10FEFS,RWYK20RMTI,TSIM12RFEM,VaRyKa15EMMI,Wick13FEFS}
shows a high computational efficiency, although the supporting rigorous
analysis  or an attempt to pass to the desired limit in the diffuse mushy
interface are  missing.
It is quite generally understood that the truly general FSI problem at large
strains is troublesome for many reasons and some regularization or linearization
is needed for a rigorous analysis, cf.\ the analysis for an incompressible
linearized fully Eulerian variant in \cite{LiuWal01EDFC}. Also, the opposite
fully Lagrangian alternative would not make  the  life easier:
although such description would allow for evolving the outer boundary of
the domain, it would need a higher gradient in the stored energy
\cite{KruRou19MMCM} which would not be conceptually consistent in the fluidic
part where the shear elastic response should completely vanish. Moreover, in such
Lagrangian setting,
the frame-indifferent viscosity would have to be very nonlinear and also the
interaction of spacial field (as gravitational) would be more complicated,
cf.\ e.g.\ the incompressible fluid interacting with elastic solid but only in
finite-dimensional approximation or small perturbation from rigid body locally
in time until topology changes in \cite{Boul07EWST,BoScTa12ESSM}. For
 all  these reasons, we will use higher gradients, involved here 
in viscosity  instead of stored energy.

The main attributes (and in their combination also novelty) of the devised
model are:\\
\Item{{\bf---}}{Concept of {\it hyperelastic solid materials}
(whose conservative-stress response comes from a stored energy) combined with
 viscosity in  the {\it Kelvin-Voigt viscoelastic rheology}.}
\Item{{\bf---}}{The fully {\it Eulerian rate formulation} of the solid part in
terms of velocity and deformation gradient is used while the deformation itself
does not explicitly occur, although it can be reconstructed a-posteriori.}
\Item{{\bf---}}{Viscoelastic  compressible  fluid formulated compatibly with
the solids but, of course, exhibiting zero shear elastic response.}
\Item{{\bf---}}{The {\it frame indifference} of the stored energy (which  has to be 
{\it nonconvex} in terms of deformation gradient) and
a physically desirable {\it singularity} under infinite compression in
relation with {\it local non-interpenetration} allowed.}
\Item{{\bf---}}{The nonconservative part of the stress in the Kelvin-Voigt model
containing a higher-order component reflecting the concept of nonsimple
{\it multipolar media} is exploited.}
\Item{{\bf---}}{The model allows for rigorous mathematical analysis as far as
 global  existence and certain regularity of weak solutions
concerns.}

To highlight the main concepts and phenomena, we confine ourselves on spatially
homogeneous fluids and homogeneous solids (except Remark~\ref{rem-inhomogenous})
and on quasistatic models (except Remark~\ref{rem-inertia})  which
is legitimate in regimes not creating elastic waves or too fast vibrations,
in particular for external loading varying only slowly in time and
for initial conditions not far from steady states. We also confine ourselves on isothermal
situations, although the energy dissipation equality as claimed in
Proposition~\ref{prop-Euler}(ii) indicates that fully thermodynamical
extension would be analytically {amenable}, too.

The plan of this article is: First, we specify the models for viscoelastic
solids and fluids and their mechanical interaction 
in Section~\ref{sec-model}. Then we reformulate this model
in a unified way (called ``{monolithic}'') in Section~\ref{sec-mono}
and show its energetics.
Weak solutions are then defined and their existence are proved in
Section~\ref{sec-anal} by the Schauder fixed point argument combined with
a suitable cut-off regularization, when using also nontrivial arguments
to cope with sharp interface between solid and fluidic domains.

\noindent
The main notation used in this paper is summarized in the following table:
\begin{center}
\fbox{
\begin{minipage}[t]{16em}\small\smallskip
$\varOmega$ a fixed domain in $\R^d$, $d=2,3$,\\
$\varGamma$ the boundary of $\varOmega$,\\
$\OmegaF,\,\OmegaS$ open subsets of $\Omega$,\\
$\xx\in\varOmega$ actual (Eulerian) coordinate,\\
$\FF$ deformation gradient (Eulerian),\\
$\TT$ Cauchy stress (symmetric - in Pa),\\
$\vv$ velocity (in m/s),\\
$\varrho$ mass density (in kg/m$^3$),\\
$\rhoR=\rhoR(\XX)$ referential mass density,\\ 
$\bm\xi$ return mapping (in m),\\
$\AA=\FF^{-1}$ distortion (Eulerian),\\
$\nu_\flat^{}>0$ a boundary viscosity,\\
$\R_{\rm sym}^{d\times d}=\{A\in\R^{d\times d};\ A^\top=A\}$,\\
$\bbI\in\R_{\rm sym}^{d\times d}$ the unit matrix,\\
${\rm Cof}(\cdot)$ cofactor matrix,
\smallskip \end{minipage}
\begin{minipage}[t]{21em}\small\smallskip
$\nn$ the unit outward normal to $\varGamma$,\\
 $\XX\in\varOmega$ referential (Lagrangian) coordinate,\\
$J=\det\FF$ Jacobian deformation gradient,\\
$p$ pressure in fluid (in J/m$^3$=Pa),\\
$\ee(\vv)=\frac12\Nabla\vv^\top\!+\frac12\Nabla\vv$ small strain rate (in s$^{-1}$),\\
$\DIS=\DIS(\ee(\vv))$ dissipative part of Cauchy stress,\\
$\mathfrak{H}(\nabla\ee(\vv))$ hyperstress (in Pa\,m),\\
$\varphi_\text{\sc r}^{}=\varphi_\text{\sc r}^{}(\XX,\FF)$ referential stored energy (in Pa),\\
$(^{_{_{\bullet}}})'$ (partial) derivative,\\
$(^{_{_{\bullet}}})\!\DT{^{}}=\pdt{}{^{_{_{\bullet}}}}+(\vv{\cdot}\Nabla)^{_{_{\bullet}}}$
convective time derivative,\\
$\cdot$ or $:$ scalar products of vectors or matrices,\\ 
$\Vdots\ \ $ scalar products of 3rd-order tensors,\\
$\nu>0$ a bulk hyper-viscosity coefficient,\\
$\GRAVITY$ external bulk load (gravity acceleration in m/s$^{2}$),
$\det(\cdot)$ determinant of a matrix.
\smallskip \end{minipage}
}\end{center}

\vspace{-1.em}

\begin{center}
{\small\sl Table\,1.\ }
{\small
Summary of the basic notation used through Sections~\ref{sec-model}--\ref{sec-anal}. 
}
\end{center}

\section{The mechanical model}\label{sec-model}

It is important to distinguish carefully the referential and
the actual time-evolving coordinates. Our aim is to formulate
the model eventually in actual configurations, i.e.\ the Eulerian
formulation, reflecting also the reality in many (or even most)
situations (and a certain general agreement) that a reference
configuration is only an artificial construction and, even
if relevant in some situations, becomes successively more and
more irrelevant during evolution at truly large
deformations (especially in fluids) and large (finite) strains.
On the other hand,  available  experimental material data 
are  typically  referen{}tial -- in particular 
it concerns the mass density and the stored energy per mass (in J/kg)
or per referential volume (in J/m$^3$=Pa), as considered here.

\subsection{Finite-strain kinematics and mass and momentum transport}
\label{sec-kinem}
We will present briefly the fundamental concepts and formulas which can
mostly be found in many monographs, as e.g.\ \cite[Part~XI]{GuFrAn10MTC} or
 \cite[Sect.~7.2]{Mart19PCM}. 
In finite-strain continuum mechanics, the basic geometrical concept is the
time-evolving deformation $\yy:\varOmega\to\R^d$ as a mapping from a reference
configuration of the body $\varOmega\subset\R^d$ into a physical space $\R^d$.
The ``Lagrangian'' space variable in the reference configuration will be
denoted as $\XX\in\varOmega$ while in the ``Eulerian'' physical-space
variable by $\xx\in\R^d$. The basic geometrical object is the (referential)
deformation gradient $\FF_\text{\!\sc r}^{}=\Nabla_{\!\XX}^{}\yy$.

We will be interested in deformations $\xx=\yy(t,\XX)$ evolving in time, which
are sometimes called ``motions''. The important quantity is the referential
 ``Lagrangian''  velocity $\vv_\text{\sc r}^{}=\pdt{}\yy$.
As already said, our formulation will be {purely}  Eulerian, i.e.\ in
deform{}ing  configurations. To this aim, one defines the inverse of a motion
$\bm\xi(t):\xx\mapsto\yy^{-1}(t,\XX)$, mostly called
{\it return} (alternatively called also a {\it reference}) {\it mapping}
or sometimes also as the inverse of a motion or inverse deformation map,
etc. Then we consider the actual (i.e.\ Eulerian) deformation gradient
$\FF=\FF_\text{\!\sc r}^{}{\circ}\,\bm\xi$ and the actual velocity
$\vv=\vv_\text{\sc r}^{}{\circ}\,\bm\xi$. The return mapping $\bm\xi$
satisfies the transport equation
\begin{align}\nonumber\\[-2.7em]
\DT{\bm\xi}=\bm0\,,
\label{transport-xi}\end{align}
where (and thorough the whole article) we use the dot-notation
$(\bm\cdot)\!\DT{^{}}:=\pdt{}(\bm\cdot)+(\vv{\cdot}\nabla)\bm\cdot$
for the {\it convective time derivative}
applied to scalars or, component-wise, to vectors or tensors.

Then the velocity gradient
$\Nabla\vv=\nabla_{\!\XX}^{}\vv\nabla_{\!\xx}^{}\XX=\DT\FF\FF^{-1}$,
where we used the chain-rule calculus  and
$\FF^{-1}=(\nabla_{\!\XX}^{}\xx)^{-1}=\nabla_{\!\xx}^{}\XX$. 
This gives the {\it transport-and-evolution equation for the
deformation gradient} as
\begin{align}\nonumber\\[-2.7em]
\DT\FF=(\nabla\vv)\FF\,.
  \label{ultimate}\end{align}
From this, we also obtain the evolution-and-transport equation for the Jacobian
$J=\det\FF$ as
\begin{align}\nonumber\\[-2.7em]\label{DT-det}
\DT J&={\rm Cof}\FF{:}\DT\FF=J\FF^{-\top}\!{:}\DT\FF
=J\bbI{:}\DT\FF\FF^{-1}=J\bbI{:}\Nabla\vv=J{\rm div}\,\vv\,,
\end{align}
where $\bbI$ denotes the unit matrix  and where ``Cof'' denotes
the cofactor matrix, i.e.\ the matrix composed from signed $(d{-}1){\times}(d{-}1)$ minors;
actually ${\rm Cof}\FF=(\det\FF)\FF^{-\top}$\!.
Sometimes, the stored energy
is expressed in {terms}  of the left Cauchy-Green tensor
$\BB=\FF\FF^\top$ which is then evolving/transported as
$\DT\BB=(\nabla\vv)^\top\BB+\BB(\nabla\vv)$.

The reference-mapping gradient, called a {\it distortion}, is
$\AA=\nabla\bm\xi=\FF^{-1}$ and together with its determinant $\det\AA=1/J$,
likewise \eq{ultimate} and \eq{DT-det}, satisfy the
transport-and-evolution equations
\begin{align}\nonumber\\[-2.7em]\label{distortion}
\DT\AA=-\AA(\nabla\vv)\ \ \ \text{ and }\ \ \
\DT{\overline{\det\AA}}\ =-(\det\AA){\rm div}\,\vv\,.
\end{align}
To highlight the dependence on the Eulerian coordinate $\xx$, one can express
$\FF=(\nabla\bm\xi)^{-1}$ by the standard algebra through cofactors as
\begin{align}\label{Euler-F}
\FF(\xx)=\frac{{\rm Cof}(\nabla\bm\xi(\xx))^\top}{\det(\nabla\bm\xi(\xx))}\,.
\end{align}

Basic ingredients of mechanical models are defined in the referential
frame, which also corresponds to the standardly available experimental data.
We consider a {\it hyperelastic material} in the {\it viscoelastic rheology}
of the Kelvin-Voigt type.
The basic ingredients are thus the referential mass density $\rhoR=\rhoR(\XX)$,
the elastic stored energy $\varphi_\text{\sc r}^{}=\varphi_\text{\sc r}^{}(\XX,\cdot)$,
and the  non-conservative  dissipative  stress 
$\DIS_\text{\sc r}^{}=\DIS_\text{\sc r}^{}(\XX,\cdot)$.
 Let us remind that the adjective ``hyperelastic``
refers to the (widely accepted) concept that the elastic response comes from
the mentioned stored energy $\varphi_\text{\sc r}^{}$. 
Wanting to formulate the model
in Eulerian frame, we use the shorthand notation which however indicates
dependence on $\bm\xi$ considered at a current time $t$:
\begin{align}\label{notation-Eulerian}
\!\!\rhoRxi( t,\xx)=\rhoR(\bm\xi( t,\xx)),\ \ \ 
\phiRxi( t,\xx,\cdot)=\varphi_\text{\sc r}^{}(\bm\xi( t,\xx),\cdot),
\ \text{ and }\ 
\DIS_\text{\sc r}^{\bm\xi}( t,\xx,\cdot)=\DIS_\text{\sc r}^{}(\bm\xi( t,\xx),\cdot).
\end{align}
This stored energy $\phiRxi=\phiRxi( t,\xx,\FF)$ (still counted
per referential volume) is to be consistently considered dependent on the
Eulerian  deformation gradient  $\FF$ from \eq{Euler-F} while
$\DIS_\text{\sc r}^{\bm\xi}=\DIS_\text{\sc r}^{\bm\xi}( t,\xx,\ee(\vv))$
is dependent on the symmetric gradient  of the Eulerian  velocity $\vv$.
The actual Eulerian stored
energy would then be $\varphi^{\bm\xi}( t,\xx,\FF)
=\phiRxi( t,\xx,\FF)/\!\det\FF$.
Then also the conservative part of the Cauchy stress
\begin{align}\nonumber\\[-2.6em]\label{Cauchy-stress}
\TT^{\bm\xi}( t,\xx,\FF)=\frac{[\phiRxi]_\FF'(t,\xx,\FF)\FF^\top(\xx)}{\det\FF(\xx)}
\end{align}
depends on Eulerian variables $\xx$ and $\FF$.  The formula
\eq{Cauchy-stress} is actually obtained by the variation in the reference
frame of the first Piola-Kirchhoff stress with respect to
$\FF_\text{\!\sc r}^{}$ when push forward to the actual Eulerian frame. 

The mass density (in kg/m$^3$) is a so-called extensive variable, and its
transport (expressing the conservation of mass) writes as the {\it continuity
equation} $\pdt{}\varrho+{\rm div}(\varrho\vv)=0$,
or, equivalently,  as  the {\it mass evolution-and-transport equation}
\begin{align}
\DT\varrho=-\varrho\,{\rm div}\,\vv\,.
\label{cont-eq+}\end{align}
{}Imposing  the initial condition  $\varrho|_{t=0}^{}=
\rhoR/\!\det\FF|_{t=0}^{}$ for \eq{cont-eq+}, 
one can  alternatively 
determine the density $\varrho$ instead of
the  differential  equation \eq{cont-eq+} from the algebraic relation
\begin{align}
\varrho=\frac{\rhoRxi}{\det\FF}=\rhoRxi\det(\nabla\bm\xi)
\,;\label{density-algebraically}\end{align}
{}recall that 
$\rhoR=\rhoR(\XX)$ is the mass density in the reference configuration.

All the equations \eq{transport-xi}--\eq{density-algebraically} are
understood for a.a.\ $\xx$ and are thus truly Eulerian. 
For the Eulerian formulation  of continuum mechanics 
see  classical textbooks as  e.g.\ \cite{GuFrAn10MTC,MarHug83MFE,Mart19PCM}.

\subsection{Viscoelastic  multipolar  solids}\label{sec-solid}

For simplicity, let us
consider spatially homogeneous material, cf.~Remark~\ref{rem-inhomogenous}
for an inhomogeneous case. Then the referential density and the referential
stored energy are $\XX$-independent and thus $\rhoRxi$ and
$\phiRxi$ do not depend on $\bm\xi$ and we can 
write simply $\rhoS$ and $\varphi_\text{\sc s}$, respectively, with the
subscript ``S'' standing for ``solid''. Actually, we consider
$\varphi_\text{\sc s}:{\rm GL}^+(d)\to\R$ with the group of invertible matrices
with positive determinant ${\rm GL}^+(d)=\{F\in\R^{d\times d};\ \det\,F>0\}$.
The further  mentioned  ingredient is a (possibly nonlinear) monotone
 non-conservative stress, i.e.\ here 
$\DIS_\text{\sc s}:\R_{\rm sym}^{d\times d}\to\R_{\rm sym}^{d\times d}$.
Here, in
addition, we consider also the nonsimple (multipolar) viscosity, i.e.\ the
non-conservative monotone, so-called {\it hyperstress}
$\mathcal{H}:\R^{d\times d\times d}\to\R^{d\times d\times d}:\bm{G}\to\nu|\bm{G}|^{s-2}\bm{G}$
with some $\nu>0$ and $s>d$; the
preposition ``hyper'' means that it contributes to the stress through its
divergence.
Neglecting inertia, this sort of solids are then governed by the quasistatic system
for $\vv$ and $\bm\xi$:
\begin{subequations}\label{solid}
\begin{align}\nonumber
&{\rm div}\Big(
\frac{\varphi_\text{\sc s}'(\FF)\FF^\top\!\!}{\det\FF}\,{+}\,\DIS_\text{\sc s}(\ee(\vv))
     \,{-}\,{\rm div}\mathcal{H}(\nabla\ee(\vv))\Big)
     +\frac{\rhoS\,\GRAVITY}{\det\FF}=0
\\&\label{solid1}\hspace*{4em}\text{ with}\ \ \ \FF=
\frac{{\rm Cof}(\nabla\bm\xi)^\top}{\det(\nabla\bm\xi)}\,,\ \ \ \mathcal{H}(\nabla\ee(\vv))=\nu|\nabla\ee(\vv)|^{s-2}\nabla\ee(\vv)\,,
\ \ \text{ and}
     \\&\label{solid2}\DT{\bm\xi}={\bm0}\,.
  \end{align}\end{subequations}
{}The dynamic variant with an inertial force will be discussed in
Remark~\ref{rem-inertia} below. 
The boundary conditions for the 4th-order equation \eq{solid1} are rather
delicate because of the higher-order viscosity and will be specified later.

 The classical simple-material models involve basic strain/stress (and
their rates), but it brings serious analytical troubles due to inevitable
geometrical and material nonlinearities at large strains, as articulated
in particular by J.M.\,Ball \cite{Ball02SOPE,Ball10PPNE}. 
The mentioned concept of the so-called
 {\it nonsimple materials}) represents a widely used generalization by
considering strain or stress (rates) gradients (or even fractional gradients
leading to integral-type truly nonlocal models). It opens a great menagerie
of options and of microstructural interpretations, cf.\ \cite{Stra23TCHG}.
It dates back to  the general concepts of Green, Rivlin, Mindlin, and
Toupin \cite{GreRiv64MCM,Mind64MSLE,Toup62EMCS}.  Mechanically, higher
gradients can suit for fitting various phenomenological
effects more properly than simple-material models, in particular flow profiles
\cite{BeBlNe92PBMV,BleGre67DF,FriGur06TBBC} or, in the dynamical variant as in
Remark~\ref{rem-inertia}, attenuation and speed dispersion in propagation of
elastic waves as discussed in \cite{Jira04NTCM,Roub19QSF}. 
More specifically, the higher gradients in strain rate, which is the option
adopted here under the name {\it multipolar materials}, leads to so-called
{\it normal dispersion}, i.e.\ the speed of wave decreases with their frequency.

 The original concept from 60ieth was later developed     
 by J.\,Ne\v cas at al.\ as multipolar fluids
\cite{Neca94TMF,NeNoSi89GSIC,NecSil91MVF,Novo92VMFP} or solids
\cite{NecRuz92GSIV,Ruzi92MPTM,Silh92MVMS}, later also
by E.\,Fried and M.\,Gurtin \cite{FriGur06TBBC}. In the context of fluid
with rigid solid interaction, one refers to \cite{FeiNec08MSRB}.

More specifically,  in \eq{solid1} 
we  have used the concept of  2nd-grade nonsimple  media,
in fluidic variant  also called {\it bipolar} (or dipolar) {\it fluids}
\cite{BelBlo14IBNV,BleGre67DF}.
 Here we used it  in a nonlinear variant with
the ``analytical'' goal to ensure $\vv(t,\cdot)$ Lipschitz continuous
 integrably in time, i.e.\ $\|\nabla\vv(t,\cdot)\|_{L^\infty(\varOmega;\R^{d\times d})}$
in $L^1(I)$. This 
avoids formation of singularities in the transported fields,
i.e.\ in $\bm\xi$, $\varrho$, $\FF$, $\AA$, and $ J$,
cf.\ also Remark~\ref{rem-reglar-transport} below. For this, we will need
$s>d$ to rely on the Sobolev embedding 
$W^{2,p}(\varOmega)\subset L^\infty(\varOmega)$. This qualification together
with the ``impenetrability'' boundary condition $\vv{\cdot}\nn=0$ seems
``nearly'' necessary, cf.\ \cite{Desj97LTEI} for $s\ge d$.

Reminding the mentioned difficulty of nonlinear viscoelastodynamics of simple
materials \cite{Ball02SOPE,Ball10PPNE} as far as mere existence of conventional
weak-solutions, let us mention results for Eulerian
nonsimple incompressible materials with linear viscosity of the 2nd-grade
\cite{BelBlo14IBNV} and of the 4th-grade \cite{NecRuz92GSIV};
the former one not yielding enough regular velocity as needed
here for eliminating singularities in the transport \eq{transport-xi} and
\eq{ultimate}. For nonlinear 2nd-grade viscosity in the compressible media,
we refer to \cite{Roub22VELS,Roub22QHPL,Roub22TVSE}. For completeness,
let us mention that there are some results \cite{BeNeRa99EUFM} for multipolar
variants involving also higher time derivatives \cite{DunFos74TSBF}.

It is physically reasonable to assume the frame indifference of the stored
energy, i.e.\ $\varphi_\text{\sc s}(F)=\varphi_\text{\sc s}(QF)$
for any $Q\in{\rm SO}(d)=\{A\in\R^{d\times d};\ A^\top\!A=AA^\top\!=\bbI,\ \det A=1\}$. In particular, this makes  the conservative part of  the
Cauchy stress  $\TT$  symmetric.

 An example for the frame-indifferent energy $\varphi_\text{\sc s}$ is the
so-called {\it neo-Hookean material} whose elastic response is governed by the
stored-energy
\begin{align}\label{neo-Hookean}
\varphi_\text{\sc s}(\FF)=\frac12K_\text{\sc e}^{}\big(\det\FF-1\big)^2\!+
\frac12G_\text{\sc e}^{}\Big(\frac{{\rm tr}(\FF\FF^\top)}{(\det\FF)^{2/d}}-d\Big)
\end{align}
for $\det\FF>0$ otherwise $\varphi(\FF)=+\infty$. Alternatively,
the second term in \eq{neo-Hookean} is considered as
$G_\text{\sc e}^{}({\rm tr}(\FF\FF^\top)/(\det\FF)^{2/d}-d-{\rm ln}(\det\FF))/2$ to comply
with \eq{Euler-ass-phi} below. In \eq{neo-Hookean}, $K_\text{\sc e}^{}$ means the
elastic bulk modulus and $G_\text{\sc e}^{}$ is the elastic shear modulus, also called
the second Lam\'e coefficient. This ansatz has been used e.g.\ in
\cite{BeFoIo22EFVM,DunRan06AFEA,RanRic10AFEM,TSIM12RFEM,WJWP20ESOM},
possibly in an incompressible variant as in Remark~\ref{rem-incompress} below.
Another popular stored energy describes a so-called 
{\it St.\,Venant-Kirchhoff material}, as used e.g.\ in
\cite{DunRan06AFEA,HroTur06MFEM,Piro18ESME,RanRic10AFEM}:
\begin{align}\label{SVK}
\varphi_\text{\sc s}(\FF)=\frac12\Big(K_\text{\sc e}^{}-\frac2dG_\text{\sc e}^{}\Big)({\rm tr}\EE)^2+G_\text{\sc e}^{}|\EE|^2\ \ \text{ with }\ \ \EE=(\FF^\top\!\FF-\bbI)/2
\end{align}
where $\EE$ is Green-Lagrange (sometimes also called Green-St.\,Venant) strain
tensor and $K_\text{\sc e}^{}$ is the bulk modulus; actually $K_\text{\sc e}^{}\!-2G_\text{\sc e}^{}/d$
is the first Lam\'e coefficient.

\subsection{Viscoelastic fluids}\label{sec-fluid}

Again, we confine to homogeneous fluids where the referential density and
the referential stored energy are $\XX$-independent and thus we can write
simply $\rhoF$ and $\varphi_\text{\sc f}$ instead of $\varrho_\text{\sc f}^{\bm\xi}$
and $\varphi_\text{\sc f}^{\bm\xi}$, respectively, with the subscript ``F''
standing for ``fluid''.

Fluids are characterized as continua whose
shear elastic response vanishes and whose volumetric response
does not bear negative pressures (or at least not big negative pressures).
Thus the only elastic response is in the volumetric part, i.e.\
$\varphi_\text{\sc f}$ depends only on the isochoric part $(\det\FF)\bbI$
of $\FF$ and completely ignores the deviatoric part $\FF\!/\!\det\FF$ and
we can write $\varphi_\text{\sc f}(\FF)=\phi_\text{\sc f}(J)$ with $J=\det\FF$
and some $\phi_\text{\sc f}:(0,+\infty)\to(0,+\infty)$,
cf.\ e.g.\ \cite[p.10]{MarHug83MFE}. As fluids do not
withstand (too much big) tension, $\phi_\text{\sc f}$ is not coercive.
Actually fluids are primarily liquids or gases (disregarding some other
materials as magma etc). From the mechanical viewpoint, the difference
between liquids and gases is that gases have zero mass density (vacuum) if
the pressure $p=-\phi_\text{\sc f}'(J)$ is zero, in contrast to liquids which
sometimes bear even a slightly negative pressure without going to vacuum.

Using the calculus $\det'(\cdot)={\rm Cof}(\cdot)$ and
the algebra ${A}^{-1}={\rm Cof}{A}^\top\!/\!\det{A}$
 for a square matrix $A$, 
the conservative part of the Cauchy  stress then reduces to 
\begin{align}
  \frac{\varphi_\text{\sc f}'(\FF)\FF^\top\!\!}{\det\FF}
 =\frac{\phi_\text{\sc f}'(\det\FF)(\det'\FF)\FF^\top\!\!}{\det\FF}=\frac{\phi_\text{\sc f}'(\det\FF)({\rm Cof}\,\FF)\FF^\top\!\!}{\det\FF}=\phi_\text{\sc f}'(J)\bbI\,.
\end{align}
Thus gases are characterized by $p\to0+$ for $J\to+\infty$.

The system \eq{solid} of $d{+}d$ equations thus reduces to the system of $d{+}1$
equations:
\begin{subequations}\label{fluid}
\begin{align}\nonumber
&{\rm div}\big(
\DIS_\text{\sc f}(\ee(\vv)){-}{\rm div}\mathcal{H}(\nabla\ee(\vv))\big)
     +\frac{\rhoF\,\GRAVITY}J=\nabla p
  \\&\label{fluid1}\hspace*{4em}\text{ with}\ \ \ 
p=-\phi_\text{\sc f}'(J)\,,\ \ \ \ \ \mathcal{H}(\nabla\ee(\vv))=
\nu|\nabla\ee(\vv)|^{s-2}\nabla\ee(\vv)\,,\ \ \text{ and }
\\&\label{fluid2}\DT J=J\,{\rm div}\,\vv\,.
  \end{align}\end{subequations}
Actually, the evolution-and-transport equation
$\DT J=J\,{\rm div}\,\vv$ in \eq{fluid2} can be replaced by the algebraic
relation $J=1/\!\det(\nabla\bm\xi)$ provided the transport equation \eq{solid2}
would be added. This is actually used in what follows.

Realizing that the actual density $\varrho=\rhoF/J$, we can express the
pressure as a function of density as
$p=-\phi_\text{\sc f}'(J)=-\phi_\text{\sc f}'(\rhoF/\varrho)$. This is
to be understood as an {\it isentropic state equation}.

\subsection{Fluid-solid interaction}

We consider a fixed bounded domain $\varOmega\subset\R^d$ covered
up to zero measure by two (not necessarily connected open sets
$\OmegaSt$ and $\OmegaFt$ depending
on time $t\in I=[0,T]$ with $T$ a fixed time horizon, containing
the solid and the fluid, respectively. The interface between
$\OmegaFt$ and $\OmegaSt$ is denoted
by $\varGamma_\text{\!\sc fs}(t)=\barOmega_\text{\sc f}(t)\cap
\barOmega_\text{\sc s}(t)$ with the bar denoting the closure.
The (fixed) boundary of $\varOmega$ will be denoted by $\varGamma$,
cf.\ Figure~\ref{fig1}.
\begin{center}
\begin{my-picture}{11}{3.7}{fig1}
\psfrag{Wl}{$\OmegaF$}
\psfrag{Ws}{$\OmegaS$}
\psfrag{Wl(t)}{$\OmegaFt$}
\psfrag{Ws(t)}{$\OmegaSt$}
\psfrag{G}{$\varGamma$}
\psfrag{Gsl}{$\varGamma_\text{\!\sc fs}$}
\psfrag{G(t)}{$\varGamma_\text{\!\sc fs}(t)$}
\psfrag{xi}{$\bm\xi(t)$}
\psfrag{x}{$\xx$}
\psfrag{X}{$\XX$}
\psfrag{t6}{$m_\text{\sc s}(0)$}
\hspace*{-4em}\includegraphics[width=37em]{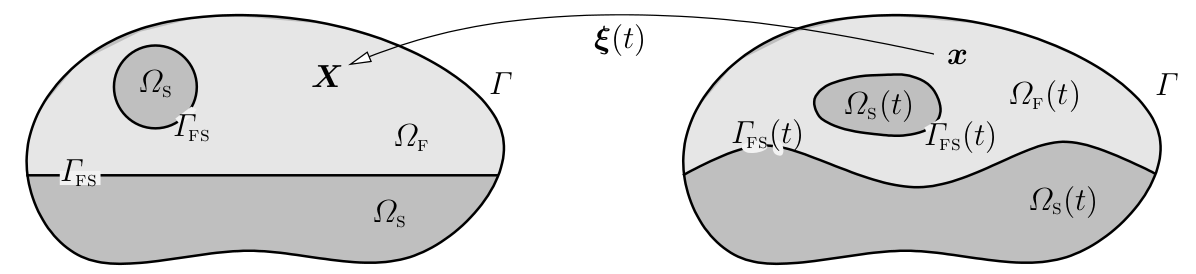}
\end{my-picture}
\nopagebreak
\\
{\small\sl\hspace*{1em}Fig.~\ref{fig1}:~\begin{minipage}[t]{33em}
A reference configuration (left) and a deformed configuration at time $t$
(right) with the boundary $\varGamma$ being fixed.
\end{minipage}
}
\end{center}

Of course,   
we need to prescribe the boundary conditions on $\varGamma$ and also the
coupling conditions on the (evolving) interface $\varGamma_\text{\!\sc fs}(t)$.
Let us note that both \eq{solid1} and \eq{fluid1} are elliptic problems for
4th-order systems for $d$-dimensional vectorial fields $\vv$ respectively on
$\OmegaSt$ and on $\OmegaFt$ parameterized by time. Therefore, they need $d{+}d$
conditions on the boundaries of each domains. On the
fluid-solid interface $\varGamma_\text{\!\sc fs}(t)$ where two domains
merge, we thus need $4d$ transition conditions.

On the boundary $\varGamma$ of $\varOmega$ which is considered fixed,
cf.\ Figure~\ref{fig1},
we denote by $\nn$ the unit outward normal to the boundary $\varGamma$
and prescribe $d{+}d$ boundary conditions: 
\begin{align}\label{BC}
\vv={\bm0}\ \ \ \text{ and }\ \ \ \Nabla\ee(\vv){:}(\nn{\otimes}\nn)={\bm0}\,.
\end{align}
The first condition fixes not only the shape of the boundary (due to the
normal velocity zero) as most frequently adopted in literature for Eulerian
formulation, but even the entire velocity field also in the tangential
direction. This is needed for global invertibility of
$\bm\xi(t,\cdot):\varOmega\to\varOmega$ used for \eq{change-of-variable} below.
The latter condition in \eq{BC} is  one of simple  
variationally consistent option{s}  for the weak formulation of the
4th-order hyper-viscosity terms in \eq{solid1} and \eq{fluid1}. 
 Another analytically well amenable option instead of Neumann-type
conditions $\Nabla\ee(\vv){:}(\nn{\otimes}\nn)={\bm0}$ would be the higher-order
Dirichlet condition $(\nn{\cdot}\nabla)\vv={\bm0}$ or a higher-order Neumann
condition $\nn{\cdot}{\rm div}(|\Nabla\ee(\vv)|^{s-2}\Nabla\ee(\vv))={\bm0}$. 
Some  options for a  relaxation of these (sometimes rather
artificial)  conditions  are  later outlined in Remarks~\ref{rem-BC+}
and \ref{rem-BC++}. 

The transition conditions on $\varGamma_\text{\!\sc fs}(t)$ are more or less
dictated by our (quite natural) intention to allow for a ``monolithic''
description of the overall fluid-solid interaction problem. More
specifically, as quite standard, we prescribe continuity of velocity
(i.e.\ {\it no-slip condition}) and {\it traction equilibrium} and,
in addition due to the multipolar material, also the 
continuity of symmetric velocity gradient, i.e.
\begin{align}\label{BC-interface}
&\Jump{\vv}=\bm{0}\,,\ \ \ \ \ \tt_\text{\sc s}=\tt_\text{\sc f}\,,\ \ \
\text{ and }\ \ \ \Jump{\ee(\vv)}=\bm{0}
\ \ \ \text{ on }\ \varGamma_\text{\!\sc fs}(t)\,,
\end{align}
where $\jump{\cdot}$ denotes the jump of normal component of the
indicated field across $\varGamma_\text{\!\sc fs}(t)$. The tractions of the
solid $\tt_\text{\sc s}$ and the fluid $\tt_\text{\sc f}$ on
$\varGamma_\text{\!\sc fs}(t)$ are respectively the vectors
\begin{subequations}\label{traction}\begin{align}\label{traction1}
&\tt_\text{\sc s}=\Big(
\frac{\varphi_\text{\sc s}'(\FF)\FF^\top\!\!}{\det\FF}\,{+}\,\DIS_\text{\sc s}(\ee(\vv))
     \,{-}\,{\rm div}\mathcal{H}(\nabla\ee(\vv))\Big)\nn_\text{\sc s}
     \ \ \text{ and}
\\&\label{traction2}
\tt_\text{\sc f}=\big(\phi_\text{\sc f}'(\det\FF)\bbI\,{+}\,
\DIS_\text{\sc f}(\ee(\vv))\,{-}\,{\rm div}\mathcal{H}(\nabla\ee(\vv))\big)
\nn_\text{\sc f}\,,
\end{align}\end{subequations}
where, at a current time $t$, we denote the ``outward'' unit
normal $\nn_\text{\sc s}(\xx)$ to $\varGamma_\text{\!\sc fs}(t)$ at $\xx\in\varGamma_\text{\!\sc fs}(t)$ oriented from the solid towards the fluid domain, and
$\nn_\text{\sc f}(\xx)=-\nn_\text{\sc s}(\xx)$.
The fields $\FF$ and $\ee(\vv)$ in \eq{traction1} consider values on $\OmegaSt$
while in \eq{traction2} they consider values on $\OmegaFt$. Actually, the
traction equilibrium $\tt_\text{\sc s}=\tt_\text{\sc f}$ means that there is
no interface tension considered in the model.

The overall fluid-solid interaction problem then reads in its classical
formulation as  
\begin{subequations}\label{FSI-}
\begin{align}\label{SFI-1}
&\text{the system }\ \eq{solid}\ \text{ on }\,\ \OmegaSt\,,
\\&\text{the system }\ \eq{fluid}\ \text{ on }\ \OmegaFt\,,
\\&\text{the transient \,conditions }\,\ \eq{BC-interface}\ \text{ on }\ \varGamma_\text{\!\sc fs}(t)\,,
\\&\text{the boundary conditions }\ \eq{BC}\ \text{ on }\ \varGamma\,.
 \end{align}\end{subequations}

\section{A monolithic model of fluid-solid interaction}\label{sec-mono}

The Eulerian formulation of the solid model and the relation
of the fluidic model with the stored energy degenerated in the
deviatoric part allows for a {\it monolithic formulation} to merge
both models and to incorporate the transient conditions
on the evolving fluid-solid interface $\varGamma_\text{\!\sc fs}(t)$.

Let us first define  a ``monolithic'' referential mass density
$\rhoR$, stored energy $\varphi_\text{\sc r}$, conservative part of the Cauchy
stress (let us denote it by $\TT_\text{\sc r}^{}$), and dissipative stress
$\DIS_\text{\sc r}$ as
\begin{subequations}\label{monolit}\begin{align}
 &\varphi_\text{\sc r}(\XX,\FF)=\begin{cases}
  \ \varphi_\text{\sc s}(\FF)&\!\!\!\text{if }\XX\!\in\!\OmegaS,\\
  \phi_\text{\sc f}(\det\FF)&\!\!\!\text{if }\XX\!\in\!\OmegaF,\end{cases}
 \hspace*{1.5em}
\TT_\text{\sc r}(\XX,\FF)=\begin{cases}
  \,\displaystyle{\frac{\varphi_\text{\sc s}'(\FF)\FF^\top\!\!}{\det\FF}}&\!\!\!\text{if }\XX\!\in\!\OmegaS,\\
\phi_\text{\sc f}'(\det\FF)\bbI&\!\!\!\text{if }\XX\!\in\!\OmegaF,\end{cases}
\\&
\DIS_\text{\sc r}(\XX,\ee)=\begin{cases}
  \DIS_\text{\sc s}(\ee)&\!\!\text{if }\XX\!\in\!\OmegaS\,,\\
  \DIS_\text{\sc f}(\ee)&\!\!\text{if }\XX\!\in\!\OmegaF\,,\end{cases}
  \ \ \ \text{ and }\ \ \ \rhoR(\XX)=\begin{cases}\rhoS&\!\!\text{if }\XX\!\in\!\OmegaS\,,\\\rhoF&\!\!\text{if }\XX\!\in\!\OmegaF\,.\end{cases}
\end{align}\end{subequations}
Actually, the reference 
(presumably ``small'') 2nd-grade hyper-viscosity coefficient $\nu$ can
also be different for the fluid and the solid.

As now $\varphi_\text{\sc r}^{}$, $\TT_\text{\sc r}^{}$, $\DIS_\text{\sc r}^{}$,
and $\rhoR$ are inhomogeneous, their Eulerian representation will depend on
$\bm\xi$. Therefore, from now on, we will use the notation
\eq{notation-Eulerian}, i.e.\ 
the quantities $\varphi_\text{\sc r}^{\bm\xi}$, $\TTRxi$, $\DIS_\text{\sc r}^{\bm\xi}$
and $\rhoRxi$ are thus defined for a.a.\ $\xx\in\varOmega$.

In terms of \eq{monolit}, for all $t\in I$, the systems \eq{solid} and
\eq{fluid} can be written ``monolithically'' as a single system of $d{+}d$
equations for $\vv$ and $\bm\xi$ of $\OmegaFt\cup\OmegaSt$:
\begin{subequations}\label{FSI}
\begin{align}\nonumber
&{\rm div}\bigg(
\TTRxi\Big(\frac{{\rm Cof}(\nabla\bm\xi)^\top}{\det(\nabla\bm\xi)}\Big)
 {+}\DISRxi(\ee(\vv)) {-}{\rm div}\mathcal{H}(\nabla\ee(\vv))\bigg)
     +\rhoRxi\,\GRAVITY\det(\nabla\bm\xi)=0
\\&\label{FSI1}\hspace{12em}\text{with}\ \ \ 
\mathcal{H}(\nabla\ee(\vv))=
\nu|\nabla\ee(\vv)|^{s-2}\nabla\ee(\vv)\ \ \ \text{ and}
     \\&\label{FSI2}
\pdt{\bm\xi}=-(\vv{\cdot}\nabla)\bm\xi\,;
  \end{align}\end{subequations}
{here} and in what follows,  the $\xx$-dependence of $\TTRxi$ and $\DISRxi$ and of $\rhoRxi$ 
is omitted for notational simplicity.
This monolithic model is to be completed by the boundary \eq{BC},
while it hides the fluid-solid interface $\varGamma_\text{\sc fs}^{}(t)$
as well as the transient conditions \eq{BC-interface}, cf.\
Remark~\ref{rem-BC} below. 

The energetics behind \eq{FSI} can be revealed by
testing \eq{FSI1} by $\vv$ and by using Green's formula twice on $\varOmega$
together with a surface Green formula on
$\varGamma$ because of the multipolar viscosity, cf.\
\eq{Euler-test-momentum++} below.

Let us first analyze the conservative part ${\rm div}\,\TTRxi(\FF)$ with
$\FF={\rm Cof}(\nabla\bm\xi)^\top\!/\!\det(\nabla\bm\xi)$. Using the
algebra $\FF^{-1}\!={\rm Cof}\,\FF^\top\!/\!\det\FF$ and the calculus
$\det'(\FF)={\rm Cof}\,\FF$, we can write the conservative part of the Cauchy
stress as
\begin{align}\nonumber
\TTRxi(\FF)&\ =\ \frac{[\phiRxi]_\FF'(\FF)}{\det\FF}\FF^\top
=\frac{[\phiRxi]_\FF'(\FF)
-\phiRxi(\FF)\FF^{-\top}\!\!\!\!}{\det\FF}\FF^\top\!+
\frac{\phiRxi(\FF)}{\det\FF}\bbI
\\&=\bigg(\frac{[\phiRxi]_\FF'(\FF)}{\det\FF}
-\frac{\phiRxi(\FF){\rm Cof}\FF}{(\det\FF)^2}\bigg)\FF^\top\!+
\frac{\phiRxi(\FF)}{\det\FF}\bbI
=\Big[\frac{\phiRxi(\FF)}{\det\FF}\Big]_{\!\FF}'\FF^\top\!
+\frac{\phiRxi(\FF)}{\det\FF}\bbI\,.
\label{referential-stress}\end{align}
Using the calculus \eq{referential-stress} and the matrix algebra
$A{:}(BC)=(B^\top\!A){:}C=(AC^\top){:}B$ for any square matrices $A$, $B$, and $C$
and also the evolution-and-transport equation \eq{ultimate}, we obtain
\begin{align}\nonumber
\TTRxi(\FF){:}\ee(\vv)&
=\frac{[\phiRxi]_\FF'(\FF)\!}{\det\FF}\FF^\top\!{:}\ee(\vv)
=\Big(\Big[\frac{\phiRxi(\FF)}{\det\FF}\Big]_{\!\FF}'\FF^\top\!
+\frac{\phiRxi(\FF)}{\det\FF}\bbI\Big){:}\ee(\vv)
\\&\nonumber
=\Big[\frac{\phiRxi(\FF)}{\det\FF}\Big]_{\!\FF}'
{:}(\Nabla\vv)\FF
+\frac{\phiRxi(\FF)}{\det\FF\!}\,{\rm div}\,\vv
\\&=\Big[\frac{\phiRxi(\FF)}{\det\FF}\Big]_{\!\FF}'
{:}\Big(\pdt\FF+(\vv{\cdot}\nabla)\FF\Big)
+\frac{\phiRxi(\FF)}{\det\FF\!}\,{\rm div}\,\vv\,.
\label{Euler-large-FSI-}\end{align}
{Let us recall that we agreed to omit 
the dependence on the variable $\xx$ for notational simplicity.
Further, we use the calculus
\begin{align}\nonumber
&\pdt{}\Big(\frac{\phiRxi(\FF)}{\det\FF}\Big)=
\frac{\big[[\varphi_\text{\sc r}]_\XX'\big]^{\bm\xi}(\FF)}{\det\FF\!}{\cdot}\pdt{\bm\xi}+
\Big[\frac{\phiRxi(\FF)}{\det\FF}\Big]_\FF'{:}\pdt\FF\ \ \ \ \text{ and}
\\&\nonumber\nabla\Big(\frac{\phiRxi(\FF)}{\det\FF}\Big){\cdot}\vv=
\frac{\big[[\varphi_\text{\sc r}]_\XX'\big]^{\bm\xi}(\FF)}{\det\FF\!}{\cdot}(\vv{\cdot}\nabla)\bm\xi+\Big[\frac{\phiRxi(\FF)}{\det\FF}\Big]_\FF'{:}(\vv{\cdot}\nabla)\FF\,,
\end{align}
so that \eq{Euler-large-FSI-} turns into
\begin{align}\nonumber
\TTRxi(\FF){:}\ee(\vv)&=
\pdt{}\Big(\frac{\phiRxi(\FF)}{\det\FF}\Big)+
\nabla\Big(\frac{\phiRxi(\FF)}{\det\FF}\Big){\cdot}\vv
\\&\qquad+
\frac{\phiRxi(\FF)}{\det\FF\!}\,{\rm div}\,\vv
-\frac{\big[[\varphi_\text{\sc r}]_\XX'\big]^{\bm\xi}(\FF)\!}{\det\FF}{\cdot}
\Big(\!\!\!\lineunder{\pdt{\bm\xi}+(\vv{\cdot}\nabla)\bm\xi}{$=0$ due to \eq{FSI2}}\!\!\!\Big)
\,.
\label{Euler-large-FSI}\end{align}
Analytically, it is rather nontrivial that the last term indeed
vanishes, cf.\ Step~4 in the proof of Proposition~\ref{prop-Euler} below.
Anyhow, continuing formally at this moment and relying further on the
integrability rigorously proved later in Section~\ref{sec-anal}, we can
integrate \eq{Euler-large-FSI} over $\varOmega$ and use it for the following
calculations (exploiting Green's formula twice on $\varOmega$) to obtain:
\begin{align}\nonumber
&\int_\varOmega\!-\,{\rm div}\,\TTRxi(\FF){\cdot}\vv\,\d\xx
=\int_\varOmega\TTRxi(\FF){:}\ee(\vv)\,\d\xx
-\int_\varGamma\vv{\cdot}\TTRxi(\FF)\nn\,\d S
\\&\nonumber\qquad
=\frac{\d}{\d t}\int_\varOmega\frac{\phiRxi(\FF)}{\det\FF}\,\d\xx+
\int_\varOmega\!\nabla\Big(\frac{\phiRxi(\FF)}{\det\FF}\Big){\cdot}\vv
+\frac{\phiRxi(\FF)}{\det\FF\!}\,{\rm div}\,\vv\,\d\xx
-\!\int_\varGamma\!\vv{\cdot}\TTRxi(\FF)\nn\,\d S
\\&\qquad
=\frac{\d}{\d t}\int_\varOmega\frac{\phiRxi(\FF)}{\det\FF}\,\d\xx+
\int_\varGamma\frac{\phiRxi(\FF)}{\det\FF}(\hspace*{-.7em}
\lineunder{\vv{\cdot}\nn}{$=0$}\hspace*{-.7em})
-\big[\TTRxi(\FF)\nn\big]_\text{\sc t}^{}{\cdot}\,\vv\,\d S\,.
\label{Euler-large-FSI+}\end{align}

The further contribution is from the dissipative part of the Cauchy stress, let
us abbreviate it as $\mathcal{D}^{\bm\xi}=\DISRxi(\ee(\vv))
-{\rm div}\,\mathcal{H}(\nabla\ee(\vv))$ with the hyperstress
$\mathcal{H}(\nabla\ee(\vv))=\nu|\nabla\ee(\vv)|^{ s-2}\nabla\ee(\vv)$.
It uses Green's formula over $\varOmega$ twice and the surface Green formula
over the boundary $\varGamma$ assumed smooth. Then
\begin{align}\nonumber
&\int_{\varOmega}-{\rm div}\mathcal{D}^{\bm\xi}{\cdot}\vv\,\d\xx
=\int_{\varOmega}\mathcal{D}^{\bm\xi}{:}\Nabla\vv\,\d\xx
-\int_{\varGamma}\vv{\cdot}\mathcal{D}^{\bm\xi}\nn\,\d S
\\&\qquad\nonumber
=\int_{\varOmega}\DISRxi(\ee(\vv)){:}\ee(\vv)+\mathcal{H}(\nabla\ee(\vv))
\Vdots\Nabla^2\vv\,\d\xx
-\int_{\varGamma}\vv{\cdot}\mathcal{D}^{\bm\xi}\nn-
\mathcal{H}(\nabla\ee(\vv))\Vdots(\nn{\otimes}\nabla\vv)\,\d S
\\&\qquad=
\int_{\varOmega}\!\DISRxi(\ee(\vv)){:}\ee(\vv)+\nu|\ee(\vv)|^s\,\d\xx
-\!\int_{\varGamma}\!\Big(\mathcal{D}^{\bm\xi}\nn
+\divG(\mathcal{H}(\nabla\ee(\vv))\nn)\Big){\cdot}\vv
\,\d S\,.
\label{Euler-test-momentum++}\end{align}
Here we used the surface Green formula and the boundary condition
$\mathcal{H(\nabla\ee(\vv))}{:}(\nn{\otimes}\nn)=\bm0$ for
\begin{align}\nonumber
&\int_{\varGamma}\mathcal{H}(\nabla\ee(\vv))\Vdots(\nn{\otimes}\nabla\vv)\,\d S=
\int_{\varGamma}\big(\hspace{-.7em}\lineunder{\mathcal{H}(\nabla\ee(\vv)){:}(\nn{\otimes}\nn)}{$=0$ due to \eq{BC}}\hspace{-.7em}\big){\cdot}\frac{\pl\vv}{\pl\nn}
+\mathcal{H}(\nabla\ee(\vv))\Vdots(\nn{\otimes}\nablaG
\vv)\,\d S
\\&\nonumber\
=\int_{\varGamma}
(\divG\nn)\big(\hspace{-.7em}\lineunder{\mathcal{H}(\nabla\ee(\vv)){:}(\nn{\otimes}\nn)}{$=0$ due to \eq{BC}}\hspace{-.7em}\big){\cdot}\vv-\divG(\mathcal{H}(\nabla\ee(\vv))\nn){\cdot}\vv\,\d S
-\!\!\int_{\varGamma}\!\divG(\mathcal{H}(\nabla\ee(\vv))\nn){\cdot}\vv\,\d S\,,
\end{align}
where we used the decomposition of $\Nabla\vv$ into its normal and tangential
parts, i.e.\ written as  $\Nabla\vv=(\pl\vv/\pl\nn)\nn+\nablaG\vv$.
Here, $\divG={\rm tr}(\nablaG)$ denotes the $(d{-}1)$-dimensional
surface divergence with ${\rm tr}(\cdot)$ being the trace of a
$(d{-}1){\times}(d{-}1)$-matrix and
$\nablaG v=\nabla v-({\partial v}/{\partial\nn})\nn$.

Now we sum \eq{Euler-large-FSI+} and \eq{Euler-test-momentum++}.
The last integral in \eq{Euler-large-FSI+} summed with the last integral in
\eq{Euler-test-momentum++} allows us to use the latter boundary condition
in \eq{BC}. Altogether, we thus obtain
(at least formally) the mechanical energy dissipation balance:
\begin{align}
&\hspace*{0em}\frac{\d}{\d t}\int_\varOmega\!\!\!\!
  \linesunder{\frac{\phiRxi(\FF)}{\det\FF}}{stored}{energy}\!\!\!\d\xx
+\!\int_\varOmega\!\!\!\!\lineunder{\DISRxi(\ee(\vv)){:}\ee(\vv)
+\nu|\Nabla\ee(\vv)|^s_{_{_{}}}\!}{dissipation rate due to viscosity}\!\!\d\xx
=\int_\varOmega\!\!\!\!\!\!\!\linesunder{\frac{\rhoRxi\,\GRAVITY{\cdot}\vv}{\det\FF}
}{power of}{gravity field}\!\!\!\!\!\!\d\xx\,.
\label{mech-engr}
\end{align}

\begin{remark}[{\sl  Dynamic problems: adding  inertia }]\upshape
\label{rem-inertia}
The dynamical  modification of the  pre\-sented  quasistatic models by
considering also inertia 
leads to expansion of the first equations in \eq{solid1} and \eq{fluid1}
respectively as
\begin{align*}&\pdt{(\varrho\vv)}=\Big(
\frac{\varphi_\text{\sc s}'(\FF)\FF^\top\!\!}{\det\FF}+\DIS_\text{\sc s}(\ee(\vv))
     -{\rm div}\mathcal{H}(\nabla\ee(\vv))-\varrho\vv{\otimes}\vv\Big)
     +\varrho\,\GRAVITY\ \ \text{ with}\ \ \varrho=\frac{\rhoS}{\det\FF}
\ \text{ and}
     \\&\pdt{(\varrho\vv)}={\rm div}\big(
 \DIS_\text{\sc f}(\ee(\vv))-{\rm div}\mathcal{H}(\nabla\ee(\vv))
 -\varrho\vv{\otimes}\vv\big)+\nabla\phi_\text{\sc f}'(J)
      +\varrho\,\GRAVITY\ \ \text{ with}\ \ \varrho=\frac{\rhoF}J\,.
\end{align*}
Naturally, this expansion needs an additional initial condition for $\vv$. 
Such expanded model allows for propagation of elastic waves. In solid regions,
there can be  both  longitudinal (pressure) waves and shear waves
(denoted respectively as P- and S-waves) while in fluidic regions,
only P-waves can propagate. On the fluid-solid interface, various reflections of
S-waves in solids and refractions and transformations of P-waves
can occur, as well as  so-called  Love, Rayleigh, or head (von Schmidt)
waves.  Adding  inertia would help in a-priori estimation strategy
\eq{est-of-rho.f.v}  working for more general $\GRAVITY$ and,  in
the modification of Remark~\ref{rem-BC+} below, 
would allow for pure slip on $\varGamma$ (i.e.\ $\nu_\flat^{}=0$) even if
$\TRACTION\ne0$.  Also, the mentioned improved estimation strategy
\eq{est-of-rho.f.v} by considering inertia would facilitate an anisothermal extension
like in \cite{Roub22TVSE}. On the other hand, handling the inertial force $\varrho\DT\vv$,
i.e.\ the additional terms $\pdt{}(\varrho\vv)+{\rm div}(\varrho\vv)$, brings
various technicalities in convergence of the momentum $\varrho\vv$ and the velocity
$\vv$ itself and, mainly, the uniqueness of the response $\vv$ for fixed
$\bm\xi$ is troublesome and thus the Schauder fixed-point arguments used
below do not work. Also e.g.\ \eq{Euler-est-formal1+} is not at disposal.
As a result,  the analysis would need to use a semi-Galerkin approximation
in the spirit of Remark~\ref{rem-construct} below  like in
\cite{Roub22TVSE,RouSte??VSPS}.
\end{remark}

\begin{remark}[{\sl Spatially inhomogeneous media}]\label{rem-inhomogenous}
\upshape
The homogeneity of particular solid and fluid areas in their reference
configuration was considered in Section~\ref{sec-solid} and \ref{sec-fluid}
rather for notational simplicity and smoothly inhomogeneous media
can be easily considered, as well. Actually,  when  admitting
a ``thin'' mushy-like layer instead of a sharp interface
$\varGamma_\text{\sc fs}^{}$ between fluid and solid regions so that
$\varphi_\text{\sc r}^{}(\cdot,\FF)\in W^{1,1}(\varOmega)$, the last term in
\eq{Euler-large-FSI} would be zero a.e.\ pointwise and we can even avoid the
requirement of global invertibility of $\bm\xi$ used in the proof of
Proposition~\ref{prop-Euler} to justify that this last term in
\eq{Euler-large-FSI} vanishes at least when integrated over $\varOmega$
in \eq{Euler-large-FSI+}. This would allow for considering a more general
boundary conditions \eq{BC} on $\varGamma$ in Remark~\ref{rem-BC+} below.
The inhomogeneous modification of the
mass density $\varrho_\text{\sc r}(\XX)$ and the 
dissipati{}ve stress  $\DIS_\text{\sc r}(\XX,\ee)$ is even simpler.
\end{remark}

\begin{remark}[{\sl Usage of distortion}]\upshape
In view of \eq{Euler-F}, we implicitly formulated the problem in terms of the
distortion $\AA=\nabla\bm\xi$ instead of $\FF=(\nabla\bm\xi)^{-1}$
like e.g. in \cite{PaPeKlHCM}. Let us realize that the continuity equation 
\eq{density-algebraically} can also be written as $\varrho=\rhoRxi\det\AA$
and the stress $\varphi_\text{\sc s}'(\FF)\FF^\top\!/\!\det\FF$ in \eq{solid1}
as $(\det\AA)\varphi_\text{\sc s}'(\AA^{-1})\AA^{-\top}\!$.
\end{remark}

\begin{remark}[{\sl Level-set method}]\upshape
In literature as e.g.\ \cite{ISTT11IFEM,LiuWal01EDFC,SuXuZh14FEFE,VaRyKa15EMMI},
one also play{}s  with the  tensorial  evolution-and-transport
\eq{ultimate} of $\FF$ or of $\BB=\FF\FF^\top$ and replace the vectorial transport
\eq{transport-xi} of $\bm\xi$ by a transport of a scalar $[-1,1]$-valued ``phase-field''
variable $\chi$, i.e.\ $\DT\chi=0$. Its positive (resp.\ negative) value 
decides whether the medium is solid (resp.\ fluid), while
$\varGamma_\text{\sc fs}(t)=\{\xx\in\varOmega;\ \chi(t,\xx)=0\}$.
Then, instead of \eq{monolit}, one can then make $\varphi$, $\TT$, and
$\DIS$ dependent on $\chi$ in a discontinuous way instead of $\xx$.
The analysis seems however complicated due to this discontinuity, unless
one admits some ``mushy'' interface by admitting the $\chi$-dependence
of $\varphi$, $\TT$, and $\DIS$ continuous.
\end{remark}

\begin{remark}[{\sl Incompressible models}]\label{rem-incompress}\upshape
Often, the compressibility is small and is thus neglected, which is legitimate
in particular in fluids. It means that $\det\FF=1$. In solids, this is a
nonlinear constraint while in fluids $J=1$ is affine. Anyhow, monolithically
this constraint can be ensured by the linear constraint ${\rm div}\,\vv=0$
provided $\det\FF|_{t=0}^{}=1$, i.e.\ $\det(\nabla\bm\xi|_{t=0}^{})=1$ on
$\varOmega$, which can be easily seen from \eq{DT-det}.
In particular, an incompressible neo-Hookean material, as considered e.g.\ in
\cite{DunRan06AFEA,DuRaRi10NSFS,HecPir17ESME,RanRic10AFEM,Wick13FEFS},
simplifies by omitting the first term in \eq{neo-Hookean}; actually
the incompressible variant is motivated in particular by the fluidic
part where the elastic shear modulus $G_\text{\sc e}^{}$ vanishes
so it suggests putting $K_\text{\sc e}^{}=+\infty$ in \eq{neo-Hookean}
as a legitimate approximation which clearly leads to $\det\FF=1$.
\end{remark}

\section{Existence of weak solutions}\label{sec-anal}

We will be interested in an initial-value problem \eq{FSI}  and provide
a proof of  global-in-time  existence and of a certain regularity of
weak solutions. To this aim, the concept of multipolar viscosity is essential
but, anyhow, still quite nontrivial and carefully ordered arguments will be needed.

We thus need to prescribe a suitable initial condition, i.e.\ here in our
quasistatic problem only for the return mapping $\bm\xi$. Notably, in terms
of \eq{ultimate}, we would
prescribe an initial condition $\FF_0$ for $\FF$ while, in terms of
\eq{distortion},  we would prescribe an initial condition $\AA_0=\FF_0^{-1}$.
For the return-mapping transport equation \eq{transport-xi}, we
prescribe the initial condition
\begin{align}\label{IC}
\bm\xi(0)=\bm\xi_0\,,
\end{align}
assuming $\nabla\bm\xi_0=\AA_0$. This implicitly determines an initial
condition for the mass density $\varrho_0=\rhoR\det(\nabla\bm\xi_0)$.

We will use the standard notation concerning the Lebesgue and the Sobolev
spaces, namely $L^p(\varOmega;\R^n)$ for Lebesgue measurable functions
$\varOmega\to\R^n$ whose Euclidean norm is integrable with $p$-power, and
$W^{k,p}(\varOmega;\R^n)$ for functions from $L^p(\varOmega;\R^n)$ whose
all derivatives up to the order $k$ have their Euclidean norm integrable with
$p$-power. We also write briefly $H^k=W^{k,2}$. The notation
$p^*$ will denote the exponent from the embedding
$W^{1,p}(\varOmega)\subset L^{p^*}(\varOmega)$, i.e.\ $p^*=dp/(d{-}p)$
for $p<d$ while $p^*\ge1$ arbitrary for $p=d$ or $p^*=+\infty$ for $p>d$.
Moreover, for a Banach space
$X$ and for $I=[0,T]$, we will use the notation $L^p(I;X)$ for the Bochner
space of Bochner measurable functions $I\to X$ whose norm is in $L^p(I)$
while $W^{1,p}(I;X)$ stands for functions $I\to X$ whose distributional
derivative is in $L^p(I;X)$. Also, $C(\cdot)$ and $C^1(\cdot)$
will denote spaces of continuous and continuously differentiable functions,
respectively. Eventually, $C_{\rm w}(I;X)$ for $X$ reflexive will denote
the space of weakly continuous functions $I\to X$.

Moreover, as usual, we will use $C$ for a generic constant which may
vary from estimate to estimate. Let us first state a result exploiting
spatial homogeneity of elastic response in the particular (solid and fluidic)
subdomains and a bounded-distortion argument, relying on that
the referential $(d{-}1)$-dimensional fluid-solid interface
$\varGamma_\text{\sc fs}^{}$ cannot ``inflate'' during the evolution to
some positive measure set:

\begin{lemma}\label{lem}
Let $\bm\xi$ be continuous and $\bm\xi(t,\cdot)\in C^1(\barOmega;\R^d)$ with
$\det(\nabla\bm\xi)>0$ on $\barOmega$. Then, for a.a.\ $(t,\xx)$, the value
$\bm\xi(t,\xx)$ belongs to the open set $\OmegaF\cup\OmegaS$.
\end{lemma}

\begin{proof}
Here we rely on a subtle analytical argument that $\bm\xi(t,\cdot)$
is of a so-called bounded distortion in the sense that
$|\nabla\bm\xi|^d/\!\det(\nabla\bm\xi)$ is bounded. Thus we have at our
disposal the well-known property \cite{HenKos14LMFD,Resh89SMBD} that pre-images
of zero-measure sets (here $\varGamma_\text{\!\sc fs}(t)$) are sets
of measure zero unless $\bm\xi(t,\cdot)$ is constant. Yet, $\bm\xi(t,\cdot)$
cannot be constant since  $\det(\nabla\bm\xi)>0$ is assumed.

As $\bm\xi$ is assumed continuous and $\OmegaF\cup\OmegaS$ is open and of a
full measure,  we can conclude  that $\bm\xi(t,\xx)\in\OmegaF\cup\OmegaS$
for a.a.\ $(t,\xx)\in I{\times}\varOmega$.
\end{proof}

Let us summarize the assumptions (with some $\varkappa>1$, $r>d$, and
$\eta>0$): 
\begin{subequations}\label{ass}
\begin{align}\nonumber
&\varOmega\ \text{ a smooth bounded domain of $\R^d$, }\ d=2,3,\
\\&\hspace{9em}\OmegaF,\OmegaS\subset\varOmega\ \text{ open},
\ \ \OmegaF\cap\OmegaS=\emptyset,\ 
\barOmega_\text{\sc f}\cup\barOmega_\text{\sc s}=\barOmega,
\\&\nonumber
\varphi_\text{\sc s}\,{\in}\,C^1({\rm GL}^+(d))\,,\ \ 
\forall F\in {\rm GL}^+(d):\ \ 
\varphi_\text{\sc s}(F)\ge\eta/(\det F)^{\varkappa-1}\ \ \text{ and}
\\&\hspace{8.8em}\label{Euler-ass-phi}
\forall Q\in{\rm SO}(d):\ \ \ \,
\varphi_\text{\sc s}(F)=\varphi_\text{\sc s}(QF)\,,
\\&
\phi_\text{\sc f}\,{\in}\,C^1((0,+\infty)),\ \ \forall J>0:\
\ \ \phi_\text{\sc f}(J)\ge\eta/J^{\varkappa-1}\!,
\\[-.0em]\nonumber
&\DIS_\text{\sc f},\DIS_\text{\sc s}\,{\in}\,C(\R_{\rm sym}^{d\times d};\R_{\rm sym}^{d\times d})\
\text{ strictly monotone and }\ \forall\ee\,{\in}\,\R_{\rm sym}^{d\times d}:\ \ 
\\&\qquad\quad
\eta|\ee|^2\le\DIS_\text{\sc f}(\ee){:}\ee\le(1{+}|\ee|^2)/\eta
\ \,\text{ and }\ \,
\eta|\ee|^2\le\DIS_\text{\sc s}(\ee){:}\ee\le(1{+}|\ee|^2)/\eta,
\label{ass-D}\\&
\GRAVITY\in L^{\infty}(I;L^{\varkappa'}(\varOmega;\R^d))\,,
\label{ass3}
\\&\bm\xi_0\in W^{2,r}(\varOmega;\R^d),\ \ \
\min_{\barOmega}\det(\nabla\bm\xi_0)>0\,,\ \ \ \bm\xi_0\big|_\varGamma^{}\text{ a
homeomorphism $\varGamma\to\varGamma$}\,.
\label{ass4}\end{align}\end{subequations}
The standard choice in literature is $\bm\xi_0(\xx)=\xx$, which
means $\FF_0=\bbI$ and $\varrho_0=\rhoR$, i.e.\ the unstretched continuum
at $t=0$, and then $\OmegaS(0)=\OmegaS$ and $\OmegaF(0)=\OmegaF$
and also $\varGamma_\text{\!\sc fs}(0)=\varGamma_\text{\!\sc fs}$. The condition
\eq{IC} covers this situation as a special case, the qualification
\eq{ass4} being then satisfied trivially.

As we consider only quasistatic problem, we do not need even
to assume positivity of the constant mass densities $\rhoR$ and $\rhoF$.
Let us note that \eq{ass4} with $r>d$ implies that
$\FF_0=(\nabla\bm\xi_0)^{-1}\in W^{1,r}(\varOmega;\R^{d\times d})$ because
\begin{align*}
\nabla\FF_0=\nabla(\nabla\bm\xi_0)^{-1}
=\nabla\frac{{\rm Cof}(\nabla\bm\xi_0)}{\det(\nabla\bm\xi_0)}=
\bigg(\frac{{\rm Cof}'(\nabla\bm\xi_0)\!}{\det(\nabla\bm\xi_0)}
-\frac{\!{\rm Cof}(\nabla\bm\xi_0)\otimes
{\rm Cof}(\nabla\bm\xi_0)}{\det(\nabla\bm\xi_0)^2}\bigg){:}\nabla^2\bm\xi_0\,,
\end{align*}
from which one can see that $\nabla\FF_0\in L^r(\varOmega;\R^{d\times d\times d})$.

We will fit the definition of weak solutions to the exponents $r$ and $s$
used in \eq{FSI1} and in \eq{ass4}:

\begin{definition}[Weak solutions to \eq{FSI}]\label{def}
A triple $(\vv,\FF,\bm\xi)$ with $\vv\in L^\infty(I;W^{2,s}(\varOmega;\R^d))$
with $\vv=\bm0$ on $I{\times}\varGamma$, and $\FF\in
L^\infty(I;W^{1,r}(\varOmega;\R^{d\times d}))\,\cap\,
 W^{1,s}(I;L^{r}(\varOmega;\R^{d\times d}))$ with $\min_{I\times\barOmega}\det\FF>0$,
and $\bm\xi\in L^\infty(I;W^{2,r}(\varOmega;\R^d))
\,\cap\, W^{1,s}(I;W^{1,r}(\varOmega;\R^d))$
is called a weak solution to the initial-boundary-value
problem for the system \eq{FSI} with \eq{monolit}
and with the boundary/transient conditions \eq{BC}--\eq{BC-interface}
and the initial condition \eq{IC} if $\pdt{}\bm\xi=-(\vv{\cdot}\nabla)\bm\xi$
holds a.e.\ in $I{\times}\varOmega$
and $\bm\xi(0)=\bm\xi_0$ a.e.\ in $\varOmega$, and if the integral identity
\begin{align}
\int_\varOmega\!\bigg(\TTRxi\Big(\frac{{\rm Cof}(\nabla\bm\xi)^\top}
{\det(\nabla\bm\xi)}\Big){+}\DISRxi(\ee(\vv))\bigg){:}\ee(\widetilde\vv)
+\nu|\nabla\ee(\vv)|^{ s-2}\nabla\ee(\vv)\Vdots\,\Nabla\ee(\widetilde\vv)\,\d\xx
\nonumber\\[-.4em]
=\!\int_\varOmega\det(\nabla\bm\xi)\rhoRxi\,\GRAVITY{\cdot}\wt\vv\,\,\d\xx
\label{Euler2-weak}
\end{align}
holds for any $\widetilde\vv$  smooth with $\widetilde\vv={\bm0}$
on $\varGamma$ and for a.a.\ time instants $t\in I$, the argument $t$ being
omitted in \eq{Euler2-weak} for notational simplicity.
\end{definition}

\begin{proposition}[Existence and regularity of weak solutions]\label{prop-Euler}
Let the assumptions \eq{ass} hold for $r>d$ and $\varkappa>2$ and  $s>d$.
Then:\\
\Item{(i)}
{there exist a weak solution $(\vv,\FF,\bm\xi)$ according
Definition~\ref{def}.}
\Item{(ii)}{Moreover, this solution complies with energetics in the sense that
the energy dissipation balance \eq{mech-engr} integrated over time interval
$[0,t]$ with the initial conditions \eq{IC} hold.}
\end{proposition}

\begin{proof}
For clarity, we will divide the proof into five steps.

\medskip\noindent{\it Step 1: formal a-priori estimates}.
Let us first make formally the a-priori  estimates which
follow from the energetics \eq{mech-engr} when one
uses the assumptions \eq{ass} for $\varkappa>1$ with some $r>d$ and $s>d$.

The only difficult term is $\varrho\,\GRAVITY{\cdot}\vv$ on the right-hand
side of \eq{mech-engr}, which can be estimated by H\"older's and Young's
inequalities  for $\varkappa>2$  as 
\begin{align}\nonumber
&\!\!\int_\varOmega\frac{\rhoRxi\,\GRAVITY{\cdot}\,\vv}{\det\FF\!}\,\d\xx
\le\Big\|\frac{\rhoRxi}{\det\FF}\Big\|_{L^{\varkappa}(\varOmega)}
\|\vv\|_{L^\infty(\varOmega;\R^d)}^{}\|\GRAVITY\|_{L^{\varkappa'}(\varOmega;\R^d)}^{}
\\&
\!\!\le C_{\nu,\eta,\varkappa}^{}\|\GRAVITY\|_{L^{\varkappa'}(\varOmega;\R^d)}^{}
\bigg(\!1{+}\Big\|\frac{\rhoRxi}{\det\FF}\Big\|_{L^{\varkappa}(\varOmega)}^{\varkappa}
\bigg)+\frac\eta{2}\|\ee(\vv)\|_{L^2(\varOmega;\R^{d\times d})^{}}^2\!
+\frac\nu2\|\Nabla\ee(\vv)\|_{L^{ s}(\varOmega;\R^{d\times d\times d})}^{ s},
\label{est-of-rho.f.v}\end{align}
{}with some $C_{\nu,\eta,\varkappa}^{}$ depending on $(\nu,\eta,\varkappa)$. In
\eq{est-of-rho.f.v}  we used the Korn inequality,~i.e.\ surely
$\|\vv\|_{L^\infty(\varOmega;\R^d)}\le C(\|\Nabla\ee(\vv)\|_{L^s(\varOmega;\R^{d\times d\times d})}^{}+
\|\ee(\vv)\|_{L^2(\varOmega;\R^{d\times d})^{}}$. More in details, one is to compose
the Sobolev embedding
$W_0^{1,s}(\varOmega;\R_{\rm sym}^{d\times d})\subset L^s(\varOmega;\R_{\rm sym}^{d\times d})$
for an inequality $\|\ee\|_{L^s(\varOmega;\R^{d\times d})}^2\le C_\text{\sc s}^{}(
\|\ee\|_{L^2(\varOmega;\R^{d\times d})}^2
+\|\Nabla\ee\|_{L^s(\varOmega;\R^{d\times d\times d})}^2)
\le C_\text{\sc s}'(1+\|\ee\|_{L^2(\varOmega;\R^{d\times d})}^2
+\|\Nabla\ee\|_{L^s(\varOmega;\R^{d\times d\times d})}^s)
$ for $s\ge2$ and, relying on Dirichlet boundary conditions, the Korn inequality
$\|\Nabla\vv\|_{L^s(\varOmega;\R^{d\times d})}\le
C_\text{\sc k}^{}\|\ee(\vv)\|_{L^s(\varOmega;\R^{d\times d})}$ with some $C_\text{\sc k}^{}$, and eventually again
the Sobolev embedding $W_0^{1,s}(\varOmega;\R^d)\subset L^\infty(\varOmega;\R^d)$.
 The  term $\|\rhoRxi/\!\det\FF\|_{L^{\varkappa}(\varOmega)}^{\varkappa}$
in \eq{est-of-rho.f.v} can thus be treated by the Gronwall inequality
relying on the blow-up assumption
$\varphi_\text{\sc r}^{}(\XX,\FF)\ge\eta/(\det\FF)^{\varkappa-1}$ 
in (\ref{ass}b,c); note that $\int_\varOmega\varphi_\text{\sc r}^{\bm\xi}(\FF)/\!\det\FF\,\d\xx
\ge\eta'\|\rhoRxi/\!\det\FF\|_{L^{\varkappa}(\varOmega)}^{}$ with
$\eta'=\eta/\inf\rhoR(\varOmega)$. Of course,
the last and the penultimate terms in \eq{est-of-rho.f.v} can be
absorbed in the left-hand side of the energy balance. From \eq{mech-engr},
exploiting the boundary condition $\vv=0$ and Korn's inequality, we thus obtain 
\begin{subequations}\label{Euler-est-formal}
\begin{align}\label{Euler-est-formal1}
&\|\vv\|_{L^s(I;W^{2,s}(\varOmega;\R^d))\,\cap\,L^2(I;H^1(\varOmega;\R^d))}^{}\le C
\ \ \ \ \text{ and }\ \ \ \ \sup_{t\in I}\int_\varOmega\frac{\phiRxi(\FF(t))}{\det\FF(t)}\,\d\xx\le C\,.
\intertext{The former estimate in \eq{Euler-est-formal1} with $s>d$ in
particular means Lipschitz continuity of $\vv(t,\cdot)$ with a time-integrable
Lipschitz constant, which further implies that the transport equations copies
regularity of the initial conditions, cf.\ \cite[ Appendix 5.3]{Roub22TVSE}. 
Here, in particular,
\eq{transport-xi} with the qualification of $\bm\xi_0\in W^{1,q}(\varOmega;\R^d)$
gives $\bm\xi\in L^\infty(I;W^{1,q}(\varOmega;\R^d))\cap 
 W^{1,s}(I;L^q(\varOmega;\R^d))$ with any $1\le q<\infty$.
Conceptually, the $L^\infty(I;W^{1,q}(\varOmega;\R^d))$-estimate is obtained by testing 
\eq{transport-xi} by ${\rm div}(|\nabla\bm\xi|^{r-2}\nabla\bm\xi)$
and using Gronwall's inequality, cf.\ \cite{Roub22TVSE} for a lot of technicalities.
For the $W^{1,s}(I;L^q(\varOmega;\R^d))$-estimate, we used
$\pdt{}\bm\xi=-(\vv{\cdot}\nabla)\bm\xi\in  L^s(I;L^q(\varOmega;\R^d))\cap
L^2(I;L^{2^*q/(2^*+q)}(\varOmega;\R^d))$. 
Similarly, from the evolution-and-transport \eq{distortion}  for  
the distortion $\AA=\nabla\bm\xi$ together with the qualification
$\AA_0=\nabla\bm\xi_0\in W^{1,r}(\varOmega;\R^{d\times d})$ gives
$\AA=\nabla\bm\xi\in L^\infty(I;W^{1,r}(\varOmega;\R^{d\times d}))\cap
 W^{1,s}(I;L^r(\varOmega;\R^{d\times d}))$. For the
$L^\infty(I;W^{1,r}(\varOmega;\R^{d\times d}))$-estimate, the 
\eq{distortion} is to be tested by ${\rm div}(|\nabla\AA|^{r-2}\nabla\AA)$
while using  Gronwall's inequality now additionally handling the
right-hand side term $-\AA(\nabla\vv)$ in \eq{distortion} relying on
the regularity of $\vv$ in \eq{Euler-est-formal1} , cf.\ again
\cite{Roub22TVSE} for a lot of technicalities. For the
$W^{1,s}(I;L^r(\varOmega;\R^{d\times d}))$-estimate,
 we used \eq{Euler-est-formal1} for a comparison to see that
$\pdt{}\AA=-\AA(\nabla\vv)-(\vv{\cdot}\nabla)\AA
\in L^{ s}(I;L^r(\varOmega;\R^{d\times d}))$. 
Altogether, relying on $\bm\xi_0\in W^{2,r}(\varOmega;\R^d)$, we can conclude
the estimate}
\label{Euler-est-formal3}
&\|\bm\xi\|_{L^\infty(I;W^{2,r}(\varOmega;\R^d))\,\cap\,
 W^{1,s}(I;W^{1,r}(\varOmega;\R^d))}^{}\le C\,.
\intertext{{}Similarly, 
from the evolution-and-transport equation \eq{DT-det} for $1/\!\det\AA$ and
from the qualification \eq{ass4} of the initial condition
$1/\!\det\AA_0=1/\!\det(\nabla\bm\xi_0)$, we also obtain the boundedness of
$1/\!\det\AA$ in $L^\infty(I;W^{1,r}(\varOmega))\cap 
 W^{1,s}(I;L^r(\varOmega))$ so
that $1/\!\det\AA\in C(I{\times}\barOmega)$ and also
$\min_{I{\times}\barOmega}\det\AA>0$.
As $r>d$, this implies also estimates for $\AA^{-1}=\FF$, specifically}
&\label{Euler-est-formal2}
\|\FF\|_{L^\infty(I;W^{1,r}(\varOmega;\R^{d\times d}))\,\cap\,
 W^{1,s}(I; L^r(\varOmega;\R^{d\times d}))}\le C\ \ \text{ with }\ \ 
\min_{I\times\barOmega}^{}\det\FF>1/C\,;
\intertext{ {}cf.\ also \cite[Lemma 3.2]{RouSte??VSPS}. 
We can eventually use the  quasistatic  equation \eq{FSI1} itself:
from the $L^\infty$-estimate in \eq{Euler-est-formal3},
we can eventually improve the former estimate in \eq{Euler-est-formal1} as}
\label{Euler-est-formal1+}
&\|\vv\|_{C_{\rm w}(I;W^{2,s}(\varOmega;\R^d))}\le C\,.
\end{align}\end{subequations}

\medskip\noindent{\it Step 2: a cut-off regularization}.
Referring to the formal estimates \eq{Euler-est-formal2}, we can choose
$\LAM>0$ so small that, for any possible sufficiently regular solution, it
holds 
\begin{align}
&\det\FF>\LAM\ \ \ \ \text{ and }\ \ \ \ |\FF|<\frac1\LAM
\ \ \text{ a.e.\ on }\ I{\times}\varOmega\,.
\label{Euler-est-formal4}
\end{align}
We regularize the stress $\TTRxi(\FF)$ in \eq{FSI1} by considering
a smooth cut-off $\pi_\LAM\in C^1(\R^{d\times d})$ as
\begin{align}&\pi_\LAM(\FF)
=\begin{cases}
\qquad\qquad1&\hspace{-8em}
\text{for $\det \FF\ge\LAM$ and $|\FF|\le1/\LAM$,}
\\
\qquad\qquad0&\hspace{-8em}\text{for $\det \FF\le\LAM/2$ or $|\FF|\ge2/\LAM$,}
\\
\displaystyle{\Big(\frac{3}{\LAM^2}\big(2\det\FF-\LAM\big)^2
-\frac{2}{\LAM^3}\big(2\det\FF-\LAM\big)^3\Big)\,\times}\!\!&
\\[.2em]
\qquad\qquad\displaystyle{\times\,\big(3(\LAM|\FF|-1)^2
-2(\LAM|\FF|-1)^3\big)}\!\!&\text{otherwise}.
\end{cases}
\label{cut-off-general}
\end{align}
Here $|\cdot|$ stands for the Frobenius norm $|\FF|=(\sum_{i,j=1}^d
F_{ij}^2)^{1/2}$ for $\FF=[F_{ij}]$, which guarantees that $\pi_\LAM$
is frame indifferent.
Furthermore, we also cut-off and regularize the singular nonlinearity
$1/\!\det(\cdot)$ which is employed in the momentum equation and the
stress $\TTRxi$ as
\begin{align}\label{cut-off-det}
&{\det}_\LAM(\FF):=
\min\Big(\!\max\Big(\!\det\FF,\frac\LAM2\,\Big),\frac2\LAM\,\Big)\ \
\text{ and }
\ \ \TTRexi(\FF)
=\frac{\big[\pi_\LAM\phiRxi\big]_\FF'(\FF)\FF^\top\!\!}{\det_\LAM(\FF)}\,.
\end{align}
Note that also $[\pi_\LAM\phiRxi](\xx,\cdot)\in C^1(\R^{d\times d})$ if
$\varphi_\text{\sc r}^{}(\XX,\cdot)\in C^1({\rm GL}^+(d))$ and that 
$[\pi_\LAM\phiRxi]_\FF'$ together with the regularized Cauchy stress $\TTRexi$
are bounded, continuous, and vanish if the argument $\FF\in\R^{d\times d}$
``substantially'' violates the constraints \eq{Euler-est-formal4},
specifically:
$$
\bigg(\det\FF\le\frac\LAM2\ \ \text{ or }\ \ |\FF|\ge\frac2\LAM\bigg)
\ \ \ \Rightarrow\ \ \ \TTRexi={\bm0}\,.
$$

We then consider a regularized momentum equation \eq{FSI1} and thus the overall
monolithic system
\begin{subequations}\label{FSIreg}\begin{align}
&\!\!{\rm div}\bigg(\!\TTRexi\Big(\,
\frac{{\rm Cof}(\nabla\bm\xi)^\top\!}{{\det}_\LAM(\nabla\bm\xi)}\Big)
{+}\DIS_\text{\sc r}^{\bm\xi}(\ee(\vv))
     {-}{\rm div}\big(\nu|\nabla\ee(\vv)|^{s-2}\nabla\ee(\vv)\big)\!\bigg)
     +{\det}_\LAM(\nabla\bm\xi)\rhoRxi\,\GRAVITY=\bm0\,,
\label{FSI1reg}\!\!
\\[-.3em]&\label{FSI2reg}
\!\!\pdt{\bm\xi}=-(\vv{\cdot}\nabla)\bm\xi\,.
\end{align}\end{subequations}
Of course, for \eq{FSI1reg} we consider the boundary conditions \eq{BC}
while for \eq{FSI2reg} we consider the initial condition \eq{IC}.

\medskip\noindent{\it Step 3: solving \eq{FSIreg} by Schauder's fixed point}.
 We organize this step for a fixed point of the mapping composed from
\begin{align}\label{FP}
B\ni\vv\stackrel{\text{by \eq{FSI2reg}}}{\mapsto}\bm\xi
\ \ \ \text{ and }\ \ \ \bm\xi
\stackrel{\text{by \eq{FSI1reg}}}{\mapsto}\vv\in B
\end{align}
with  a (sufficiently large) convex set $B$ in $C_{\rm w}(I;W^{2,s}(\varOmega;\R^d))$
with $\vv=\bm0$ on $I{\times}\varGamma$. More specifically, we put 
\begin{align}\nonumber
B:=\bigg\{&\vv\in C_{\rm w}(I;W^{2,s}(\varOmega;\R^d));\ \forall_{\rm a.a.}t\in I:\ \
\vv(t)=\bm0\ \text{ on }\varGamma\ \ \text{ and}
\\[-.3em]&\nonumber\quad\|\nabla\ee(\vv(t))\|_{L^s(\varOmega;\R^{d\times d\times d})}^s
+\|\ee(\vv(t))\|_{L^2(\varOmega;\R^{d\times d})}^2
\\&\qquad\ \le
\max\Big(\frac{s'}\nu,\frac2\eta\Big)\Big(\frac{2{\rm meas}(\varOmega)}\eta
L_\LAM^2+\frac{N^{s'}\!\max(\varrho_\text{\sc s}^{s'},\varrho_\text{\sc f}^{s'})\!}
{s'\nu^{1/(s-1)}}\:\|\GRAVITY\|_{L^\infty(I;L^1(\varOmega;\R^d))}^{s'}\!\Big)\bigg\}\,,
\label{radius-}\end{align}
where $\eta$ is from \eq{ass-D}, $N$ denotes the norm of the embedding of 
$\{\vv\in W^{2,s}(\varOmega;\R^d);\ \vv|_\varGamma=\bm0\}$
normed by $\|\nabla\ee(\cdot)\|_{L^s(\varOmega;\R^{d\times d\times d})}$
into $L^\infty(\varOmega;\R^d)$, and  
\begin{align}
L_\LAM:=\sup_{\XX\in\varOmega,\,\FF\in\R^{d\times d}}\bigg|\frac{\big[\pi_\LAM\varphi_\text{\sc r}^{}(\XX,\cdot)\big]'(\FF)\FF^\top\!\!}{\det_\LAM(\FF)}\:\bigg|<+\infty\,.
\label{radius}\end{align}
The motivation of this choice of $B$ will be seen later in the estimate \eq{radius+}.

First, let us consider  a fixed  $\vv$ from $B$. Then the transport
equation \eq{FSI2reg} with $\bm\xi(0)=\bm\xi_0\in W^{2,r}(\varOmega;\R^d)$ has a
solution $\bm\xi$ in a bounded subset of $L^\infty(I;W^{2,r}(\varOmega;\R^d))\,\cap\,
 W^{1,s}(I;L^r(\varOmega;\R^d))$,
cf.\ \eq{Euler-est-formal3} with some constant $C$ depending
on  the already chosen ball  $B$  but not on the
particular choice of $\vv\in B$, and the equation \eq{FSI2reg} holds a.e.\ on
$I{\times}\varOmega$ and $\bm\xi\in C(I{\times}\barOmega;\R^d)$;
cf.\ again \cite[{}Appendix 5.3]{Roub22TVSE} as used already
for the argumentation \eq{Euler-est-formal}.
 For the mentioned uniqueness of $\bm\xi$, let us note that the
transport equation \eq{FSI2reg} is linear for a fixed $\vv$; more specifically,
considering two solutions $\bm\xi_1$ and $\bm\xi_2$
satisfying $\pdt{}(\bm\xi_1{-}\bm\xi_2)=(\vv{\cdot}\nabla)\bm\xi_2{-}(\vv{\cdot}\nabla)\bm\xi_1$
gives, when tested by $2(\bm\xi_1{-}\bm\xi_2)$, the following estimate for
a.a.\ $t\in I$:
\begin{align}\nonumber
\frac{\d}{\d t}\int_\varOmega|\bm\xi_1{-}\bm\xi_2|^2\,\d\xx&=2\!\int_\varOmega
\!\!\big((\vv{\cdot}\nabla)\bm\xi_2{-}(\vv{\cdot}\nabla)\bm\xi_1\big)
{\cdot}(\bm\xi_1{-}\bm\xi_2)\,\d\xx
\\&=\!\int_\varOmega({\rm div}\,\vv)|\bm\xi_1{-}\bm\xi_2|^2\,\d\xx
\le\|{\rm div}\,\vv\|_{L^\infty(\varOmega)}^{}\|\bm\xi_1{-}\bm\xi_2\|_{L^2(\varOmega;\R^d)}^2\,,
\end{align}
from which we conclude $\bm\xi_1=\bm\xi_2$ by Gronwall's inequality when
taking $\bm\xi_1(0)=\bm\xi_0=\bm\xi_2(0)$ into account.

Further, we use this
also for the nonhomogeneous evolution-and-transport equations for
$\AA=\nabla\bm\xi$ and for $\FF=(\nabla\bm\xi)^{-1}$ to obtain $\AA$ and $\FF$
a-priori bounded in $L^\infty(I;W^{1,r}(\varOmega;\R^{d\times d}))\cap
 W^{1,s}(I;L^r(\varOmega;\R^{d\times d}))$, cf.\ 
\eq{Euler-est-formal}b,c)  with some
constant $C$ again depending on that bounded set $B$ of $\vv$'s.
In particular, we can also see that $\AA,\FF\in C(I{\times}\barOmega;\R^d)$
and  that  the equations \eq{ultimate} and \eq{distortion} are satisfied a.e.\
on $I{\times}\varOmega$.

 Moreover,  in view of these a-priori bounds and  mentioned  uniqueness
of the response $\bm\xi$ for a current $\vv$,
it is easy to see that the dependence of $\bm\xi$
on $\vv$ is (weak,weak*)-continuous as a mapping
$L^1(I;W^{2,s}(\varOmega;\R^d))\to L^\infty(I;W^{2,r}(\varOmega;\R^d))$.
 Actually, for any weakly convergent sequence of $\vv$'s,
the corresponding weakly* sequence of $\bm\xi$'s converges, for a moment in terms of
subsequences in $L^\infty(I;W^{2,r}(\varOmega;\R^d))\cap W^{1,s}(I;W^{1,r}(\varOmega;\R^d))$
and hence strongly in $C(I{\times}\varOmega;\R^d)$, which makes it easy to
pass to the limit in the transport equation \eq{FSI2}. From the mentioned
uniqueness for $\bm\xi$, we can see that eventually the whole sequence of $\bm\xi$'s
converges,   cf.\ again \cite{Roub22TVSE}.

 Now we go on to the latter mapping in \eq{FP}, considering a fixed $\bm\xi$.
 The quasistatic equation  \eq{FSI1reg} together with the mentioned
boundary conditions represents, at each time instant $t$ and for $\bm\xi=\bm\xi(t)$
fixed, a static boundary-value problem for the 4th-order quasilinear equation
\begin{align}\nonumber
{\rm div}^2\big(\nu|\nabla\ee(\vv)|^{s-2}\nabla\ee(\vv)\big)
-{\rm div}\DIS_\text{\sc r}^{\bm\xi}(\ee(\vv))&=\ff_\LAM(\bm\xi,\nabla\bm\xi)
\\\text{ with}\ \ \ \ff_\LAM(\bm\xi,\AA)\,&\!:={\rm div}\bigg(\!\TTRexi\Big(\,
\frac{{\rm Cof}\AA^\top\!}{{\det}_\LAM(\AA)}\Big)\bigg)
+{\det}_\LAM(\AA)\rhoRxi\,\GRAVITY\,,
\label{4th-order-BVP}\end{align}
which has a solution $\vv=\vv(t)\in W^{2,s}(\varOmega;\R^d)$.  If
$\DIS_\text{\sc r}^{}(\XX,\cdot):\R_{\rm sym}^{d\times d}\to\R_{\rm sym}^{d\times d}$
would have possessed some (dissipation) potential as often assumed in literature,
the solution to \eq{4th-order-BVP} could
be proved simply by the direct minimization method, while in our general
case we can use the classical Browder-Minty theorem about surjectivity for a coercive
radially-continuous monotone operator, here  
$\vv\mapsto {\rm div}^2(\nu|\nabla\ee(\vv)|^{s-2}\nabla\ee(\vv))
-{\rm div}\DIS_\text{\sc r}^{\bm\xi}(\ee(\vv))$ with the respective boundary
conditions $\vv=\bm0$
and $\nabla\ee(\vv){:}(\nn{\otimes}\nn)=\bm0$ on $\varGamma$.
Moreover, due to the strong monotonicity of this quasilinear operator,
 the mentioned solution
is unique and  the mapping $\bm\xi\mapsto\vv$ is bounded as
$W^{1,1}(\varOmega;\R^d)\to W^{2,s}(\varOmega;\R^d)$; here it is trivial
due to the $\LAM$-cut-off regularization  \eq{cut-off-det} so that
the right-hand side $\ff_\LAM(\bm\xi,\nabla\bm\xi)$ of \eq{4th-order-BVP}
is a-priori bounded in $W^{2,s}(\varOmega;\R^d)^*$. It is important to quantify
this boundedness: testing \eq{4th-order-BVP} by $\vv$, at a current time instant
$t\in I$, we obtain the estimate 
\begin{align}\nonumber
&\nu\|\nabla\ee(\vv)\|_{L^s(\varOmega;\R^{d\times d\times d})}^s
+\eta\|\ee(\vv)\|_{L^2(\varOmega;\R^{d\times d})}^2
\le\int_\varOmega\ff_\LAM(\bm\xi,\nabla\bm\xi){\cdot}\vv\,\d\xx
\\&\nonumber\qquad
=\int_\varOmega{\det}_\LAM(\nabla\bm\xi)\rhoRxi\,\GRAVITY{\cdot}\vv
-\TTRexi\Big(\,\frac{{\rm Cof}(\nabla\bm\xi)^\top\!}{{\det}_\LAM(\nabla\bm\xi)}\Big)
{:}\ee(\vv)\,\d\xx
\le\frac{N^{s'}\!\max(\varrho_\text{\sc s}^{s'},\varrho_\text{\sc f}^{s'})\!}{s'\nu^{1/(s-1)}}\:\|\GRAVITY\|_{L^1(\varOmega;\R^d)}^{s'}\!
\\&\hspace{8em}
+\frac{\nu}s\|\nabla\ee(\vv)\|_{L^s(\varOmega;\R^{d\times d\times d})}^s\!+
\frac{2{\rm meas}(\varOmega)}\eta L_\LAM^2+\frac\eta2\|\ee(\vv)\|_{L^2(\varOmega;\R^{d\times d})}^2
\label{radius+}
\end{align}
with $L_\LAM$ from \eq{radius} and with $N$ as in \eq{radius-}.
Thus $t\mapsto\vv(t)$ belongs to $B$ defined in \eq{radius-}.

Next, we are to show that this mapping $\bm\xi\mapsto\vv$ is even continuous
while taking into account that the Nemytski\u\i\ operators in
\eq{FSI1reg} involve a discontinuity across $\varGamma_\text{\!\sc fs}^{}$ and
that we do not have any estimates on $\pdt{}\vv$ at disposal.
One ingredient is surely the mentioned strong monotonicity of the underlying
quasilinear operator. Further, we use the  Arzel\`a-Ascoli-type  compact
embedding $C_{\rm w}(I;W^{2,r}(\varOmega;\R^d))\,\cap\,
 W^{1,s}(I;L^r(\varOmega;\R^d))
\subset C(I;C^1(\barOmega;\R^d))$ 
due to the compact embedding $W^{2,r}(\varOmega)\subset C^1(\barOmega)$ for 
$r>d$; cf.\ e.g.\ \cite[Lemma~7.10]{Roub13NPDE}.  We now need to show
continuity of the mappings $\bm\xi\mapsto\det_\LAM(\nabla\bm\xi)\rhoRxi\,\bm g:
C(I;C^1(\barOmega;\R^d))\to L^1(I{\times}\varOmega;\R^d)$ and
$\bm\xi\mapsto\TTRexi({\rm Cof}(\nabla\bm\xi)^\top\!/\det_\LAM(\nabla\bm\xi))
:C(I;C^1(\barOmega;\R^d))\to L^1(I{\times}\varOmega;\R^{d\times d}))$.
This is however a bit delicate issue because $\rhoR$ and
$\TT_\text{\sc r}(\cdot,\FF)$  as well as $\DIS_\text{\sc r}(\cdot,\ee)$ 
are discontinuous across the surface $\varGamma_\text{\!\sc fs}^{}$ inside the 
 reference domain  $\varOmega$, cf.\ \eq{monolit}, so that the
usual arguments of continuity of Nemytski\u\i\ operators cannot be used.
Instead, we can use Lemma~\ref{lem}. Realizing that, having a sequence
$\{\bm\xi_k\}_{k\in\N}^{}$ converging in $C(I{\times}\barOmega;\R^d)$
to $\bm\xi$, due to Lemma~\ref{lem}, we can see that $\bm\xi(t,\xx)$ is
valued (together with an neighbourhood) in the open set $\OmegaF\cup\OmegaS$
for a.a.\ $(t,\xx)$ so that, $\rhoR(\bm\xi_k(t,\xx))\to\rhoR(\bm\xi(t,\xx))$
for a.a.\ $(t,\xx)$. Also $\nabla\bm\xi_k$ converges in
$C(I{\times}\barOmega;\R^{d\times d})$ to $\nabla\bm\xi$;
here we again use the Arzel\`a-Ascoli-type compact
embedding of $C_{\rm w}(I;W^{1,r}(\varOmega;\R^{d\times d}))\,\cap\,
 W^{1,s}(I;L^r(\varOmega;\R^{d\times d}))$ into $C(I{\times}\barOmega;\R^{d\times d})$,
 cf.\ \eq{Euler-est-formal3}. 
Thus $\det_\LAM(\nabla\bm\xi_k)\to\det_\LAM(\nabla\bm\xi)$
in $C(I{\times}\barOmega)$ since $\det_\LAM(\cdot):\R^{d\times d}\to\R$
is continuous. Altogether, $\det_\LAM(\nabla\bm\xi_k)\rhoRxik\,\bm g
\to\det_\LAM(\nabla\bm\xi)\rhoRxi\,\bm g$ a.e., and thus also in
$L^1(I{\times}\varOmega;\R^d)$ by the Lebesgue dominated-convergence
theorem. Also $\TTRexik({\rm Cof}(\nabla\bm\xi_k)^\top\!/\det_\LAM(\nabla\bm\xi_k))\to\TTRexi({\rm Cof}(\nabla\bm\xi)^\top\!/\det_\LAM(\nabla\bm\xi))$
in $L^1(I{\times}\varOmega;\R^{d\times d})$ by similar arguments.
The discontinuity of $\DIS_\text{\sc r}(\cdot,\ee(\vv))$ occurring in
the weak formulation \eq{Euler2-weak} integrated over the time interval
$I$ can be treated by the same manner.

Then we can apply the Schauder fixed-point argument for the composed mapping
$\vv\mapsto\bm\xi$ and $\bm\xi\mapsto\vv$, viz \eq{FP}. The former mapping
is considered on  the 
convex bounded closed set  $B$  from \eq{radius-}  in the Banach space
$C_{\rm w}(I;W^{2,s}(\varOmega;\R^d))$ endowed with the weak topology from
$L^s(I;W^{2,s}(\varOmega;\R^d))$ which makes  $B$  weakly compact.
As to the mapping $\bm\xi\mapsto\vv$, it is important that, due to our
cut-off $\LAM$-regularization, the mentioned a-priori bounds for $\vv$ are
entirely independent of $\bm\xi$ and hold even for $\bm\xi$'s with
$\det(\nabla\bm\xi)$ not positive. Therefore, one does not need to consider
the nonconvex constraint $\det(\nabla\bm\xi)\ge\LAM$ at this stage and can work
with the whole mentioned ball $B$ which is, of course, convex. This
gives a uniquely defined velocity field $\vv\in C(I;W^{2,s}(\varOmega;\R^d))$
with $\vv{\cdot}\nn=0$ on $I{\times}\varGamma$.
We then choose $B$ so big that the mentioned bounded set of $\vv$'s
is contained in it.
The mentioned composed mapping is single-valued and (weak,weak)-continuous,
and thus it has a fixed point $\bm\xi\in B$. Together with the corresponding
$\vv$ we thus get a weak solution to the regularized system \eq{FSIreg} with
the initial and boundary conditions \eq{BC} and \eq{IC}.

\medskip\noindent{\it Step 4: energetics rigorously}.
The estimates in Step~1  relied on the energy dissipation balance which was
derived in Section~\ref{sec-mono} rather formally. We can prove it
for \eq{FSIreg} rigorously.

Due to the cut-off regularization of \eq{FSI1reg}, we can be sure that
$\vv\in L^\infty(I;W^{2,s}(\Omega;\R^d))$. Thus, the transport of $\bm\xi$ by
\eq{FSI2reg} is done through a Lipschitz velocity field, and thus
surely $\det(\nabla\bm\xi)>0$ due to
the qualification of the initial condition \eq{ass4}. Therefore,
$\FF={\rm Cof}(\nabla\bm\xi)^\top\!/\det(\nabla\bm\xi)$ is well defined
and $\TTRexi(\FF)=[\pi_\LAM\phiRxi]_\FF'(\FF)\FF^\top\!/\!\det_\LAM(\FF)=
[\pi_\LAM\phiRxi]_\FF'(\FF)\FF^\top\!/\!\det(\FF)$ due to \eq{cut-off-det}
because, if $\det(\nabla\FF)\le\LAM/2$, then $\pi_\LAM(\FF)=0$
due to \eq{cut-off-general} and thus also $[\pi_\LAM\phiRxi]_\FF'(\FF)=0$.
Also \eq{ultimate} is at our disposal. Therefore, we can use the calculus
\eq{referential-stress}--\eq{Euler-large-FSI-}. Yet, \eq{Euler-large-FSI+}
is more delicate because it is not trivial to show that the last term in
\eq{Euler-large-FSI}, when integrated over $\varOmega$, is indeed zero -- here
we emphasize that $\varphi_\text{\sc r}(\cdot,\FF)$ is discontinuous across
$\varGamma_\text{\sc fs}$ so that $[[\varphi_\text{\sc r}]_\XX']^{\bm\xi(t)}(\FF)$
is a measure supported on $\varGamma_\text{\sc fs}(t)$, in general.

At each time instant $t\in I$, with omitting the argument $t$ for notational
simplicity and by the change-of-variable formula for
$\xx=\bm\xi^{-1}(\XX)$ and by the Green formula, we have
\begin{align}\nonumber
&\int_\varOmega\!
\frac{\big[[\varphi_\text{\sc r}]_\XX'\big]^{\bm\xi(\xx)}(\FF(\xx))\!}{\det\FF(\xx)}{\cdot}
\Big(\!\!\!\!\lineunder{\pdt{\bm\xi}+(\vv{\cdot}\nabla)\bm\xi}
{\ \ \ \ \ \ \ \ $=:\bm r(\xx)$}\!\!\!\!\Big)\,\d\xx
=\!\int_\varOmega
[\varphi_\text{\sc r}]_\XX'(\XX,\FF_\text{\!\sc r}^{}(\XX)){\cdot}
\rr(\bm\xi^{-1}(\XX))\,\d\XX
\\[-.5em]&\nonumber\hspace{9em}
=\int_\varGamma\!\varphi_\text{\sc r}(\XX,\FF_\text{\!\sc r}^{}(\XX))\,
\big(\nn{\cdot}\rr(\bm\xi^{-1}(\XX))\,\d S
\\[-.5em]&\hspace{11em}-\!\int_\varOmega\!\varphi_\text{\sc r}(\XX,\FF_\text{\!\sc r}^{}(\XX))\,
{\rm div}\big(\rr(\bm\xi^{-1}(\XX))\big)\,\d\XX=0\,,
\label{change-of-variable}
\end{align}
where we used the notation $\FF=\FF_\text{\!\sc r}^{}{\circ}\,\bm\xi$
from Section~\ref{sec-kinem}. For the  change of variable, we needed however 
the global invertibility of $\bm\xi$ for which the classical result of
J.M.\,Ball \cite[Thm.1(ii)]{Ball81GISF} is exploited; here
the local invertibility $\det(\nabla \bm\xi(t))>0$ in $\varOmega$
 together with the invertibility of $\bm\xi(t)$ on $\varGamma$ were used.
The latter property is ensured by the boundary invertibility of
 $\bm\xi_0$ assumed in \eq{ass4} together with the boundary condition
 $\vv=0$ in \eq{BC} so that $\bm\xi|_\varGamma=\bm\xi_0|_\varGamma$ stays
 constant during the evolution. Further, we needed that the residuum $\bm r$
 vanishes on $\varGamma$ and its divergence vanishes in $\varOmega$. This is
 indeed guaranteed by the transport equation \eq{FSI2reg} together with
 the mentioned regularity of its solution
 $\bm\xi\in L^\infty(I;W^{2,r}(\varOmega;\R^d))\,\cap\,
  W^{1,s}(I;L^r(\varOmega;\R^d))$ so that $\rr(t)$ has a well defined (zero)
 trace on $\varGamma$ and
\begin{align}\nonumber
{\rm div}\,\rr&=\pdt{}{\rm div}\,\bm\xi
+{\rm div}((\vv{\cdot}\nabla)\bm\xi)
=\pdt{}{\rm div}\,\bm\xi
+(\nabla\vv)^\top{:}\nabla\bm\xi+\vv{\cdot}\nabla({\rm div}\,\bm\xi)
\\&=\pdt{}({\rm tr}\AA)+(\nabla\vv)^\top{:}\AA+\vv{\cdot}\nabla({\rm tr}\AA)
\in L^{ s}(I;L^r(\varOmega))
\end{align}
due to the already mentioned regularity of the distortion 
$\AA=\nabla\bm\xi\in L^\infty(I;W^{1,r}(\varOmega;\R^{d\times d}))\cap
 W^{1,s}(I;L^r(\varOmega;\R^{d\times d}))$.

Due to the mentioned quality of $\vv$, we can legitimately test
also the dissipation part of the momentum equation in its weak formulation
\eq{Euler2-weak} by $\wt\vv=\vv\in L^\infty(I;W^{2,s}(\varOmega;\R^d))$. Thus we
can rigorously obtain the energy dissipation balance \eq{mech-engr} modified
for the regularized system as
\begin{align}
  &\hspace*{0em}\frac{\d}{\d t}
  \int_\varOmega\!\!
 \frac{\pi_\LAM(\FF)\phiRxi(\FF)\!}{\det_\LAM\FF}\,\d\xx
+\!\int_\varOmega\!\DISRxi(\ee(\vv)){:}\ee(\vv)
+\nu|\Nabla\ee(\vv)|^s\,\d\xx
=\int_\varOmega\frac{\rhoRxi\,\GRAVITY\,{\cdot}\,\vv}{\det_\LAM\FF}\,\d\xx\,.
\label{mech-engr++}
\end{align}

\medskip\noindent{\it Step 5: the original problem}.
Let us remind that $\FF={\rm Cof}(\nabla\bm\xi)^\top\!/\!\det_\LAM(\nabla\bm\xi)$
resulted as the fixed point in Step~3 lives in
$L^\infty(I;W^{1,r}(\varOmega;\R^{d\times d}))\,\cap\,
 W^{1,s}(I;L^r(\varOmega;\R^{d\times d}))$
and this space is embedded into $C(I{\times}\barOmega;\R^{d\times d})$ since
$r>d$. Similarly, $\det_\LAM\FF=1/\!\det_\LAM(\nabla\bm\xi)$ lives in
$L^\infty(I;W^{1,r}(\varOmega))\,\cap\,
 W^{1,s}(I;L^r(\varOmega))$. Therefore $\FF$
together with $\det\FF$ evolve continuously in time, being valued respectively
in $C(\barOmega;\R^{d\times d})$ and $C(\barOmega)$.
Let us recall that the initial condition $\FF_0$ complies with the bounds
\eq{Euler-est-formal4} and we used this $\FF_0$
also for the $\LAM$-regularized system.
Therefore $\FF$ satisfies these bounds not only at $t=0$ but also at least
for small times $t>0$. Yet, in view of the choice \eq{Euler-est-formal4}
of $\LAM$, this means that the $\LAM$-regularization is nonactive
and  $(\vv,\FF,\bm\xi)$ solves, at least for a small time, the
original nonregularized problem
\eq{FSI} with the initial/boundary conditions \eq{BC} and \eq{IC}.
For this solution, the a-priori $L^\infty$-bounds \eq{Euler-est-formal4} hold.
Eventually, by the continuation argument, we may see that the
$\LAM$-regularization remains inactive within the whole evolution of
$(\vv,\FF,\bm\xi)$ on the whole time interval $I$.
\end{proof}

\begin{remark}[{\sl Deformation $\yy$}]\label{rem-deformation}\upshape
The global invertibility of $\bm\xi(t)$ used in the above proof ensures also
existence a deformation $\yy$ such that $\FF=\nabla\yy$, which was the
original motivation of the model when defining the deformation gradient
tensor $\FF$.
\end{remark}

\begin{remark}[{\sl Transient conditions across $\varGamma_\text{\sc fs}(t)$}]
\label{rem-BC}\upshape
The conditions on the evolving fluid-solid interface are rather implicit
in the monolithic formulation. Clearly, $\vv(t)\in W^{1,s}(\varOmega;\R^d)
\subset C(\barOmega;\R^d)$ implies $\jump{\vv(t)}=\bm0$ and
$\ee(\vv(t))\in W^{1,s}(\varOmega;\R_{\rm sym}^{d\times d})$ for $s>d$ 
implies $\jump{\ee(\vv(t))}=\bm0$. The other condition comes from the
momentum equilibrium, when assuming the weak solution to be sufficiently
smooth. Indeed, when testing the monolithic momentum equation \eq{FSI1}
by $\wt\vv$ with a compact support in $\OmegaFt$, by the Green formula
we obtain $\int_{\OmegaFt}\varrho\,\GRAVITY{\cdot}\wt\vv
-\big(\TT_\text{\sc f}(\FF){+}\DIS_\text{\sc f}(\ee(\vv))
     {-}{\rm div}\mathcal{H}(\nabla\ee(\vv))\big){:}\ee(\wt\vv)\,\d\xx=0$.
Similarly it holds for $\OmegaSt$. Then taking $\wt\vv$ with a compact support
in $\varOmega$ and integrating \eq{FSI1} and using Green formula separately
over $\OmegaFt$ and $\OmegaSt$, we obtain the identity 
\begin{align*}
&\int_{\varGamma_\text{\sc fs}(t)}\!\!\tt_\text{\sc f}{\cdot}\wt\vv+
\tt_\text{\sc s}{\cdot}\wt\vv\,\d S=0\,,
\end{align*}
cf.\ \eq{traction}.
As $\wt\vv$ takes arbitrary values on $\varGamma_\text{\sc fs}(t)$,
it must hold $\tt_\text{\sc f}=\tt_\text{\sc s}$, as claimed in 
\eq{BC-interface}. Altogether, in the physically relevant case $d=3$,
we obtained the expected number $3+6+3=4d$ of the transient conditions
on $\varGamma_\text{\sc fs}(t)$.
\end{remark}

\begin{remark}[{\sl A constructive proof}]\label{rem-construct}\upshape
 Actually,  we can use the Galerkin approximation separately of
both the regularized momentum equation \eq{FSI1reg} and the transport
equation regularized by an $r$-Laplacian, and then make a successive
limit passage first in the transport equation \eq{FSI2} and then in
the momentum equation as in \cite{Roub22QHPL}.  Such a space discretization leading
to an initial-value problem for a system of ordinary differential equations
would give a conceptually implementable algorithm.  On the other hand, 
the proof of
Proposition~\ref{prop-Euler} would then be much more technical, which
is why the shorter but non-constructive proof based on Schauder's theorem
was here presented. 
\end{remark}

\begin{remark}[{\sl Navier-type boundary conditions:~local-in-time solutions}]
\label{rem-BC+}\upshape
The~normal~deformation $\yy{\cdot}\nn$ is usually fixed by prescribing
the normal velocity $\vv{\cdot}\nn=0$ on $\varGamma$ to avoid serious troubles
in Eulerian approach. Yet, one can think at least about making the tangential
velocity not vanishing. Zero tangential velocity was actually needed only for
the global invertibility of $\bm\xi$ used in \eq{change-of-variable}.
If one has weaker (but in the FSI-context quite usual) ambitions to prove
only local-in-time solution existence, one can admit non-zero tangential
velocity  and  rely on
that the  return mapping  $\bm\xi(t,\cdot):\varOmega\to\varOmega$ is
originally globally invertible for $t=0$ and, evolving continuously, it
stays invertible at least for sufficiently
small $t>0$. Then one can think about more general boundary conditions
not fixing the tangential velocity to zero, as illustrated in
Figure~\ref{fig2}.
\begin{center}
\begin{my-picture}{11}{3.7}{fig2}
\psfrag{Wl}{$\OmegaF$}
\psfrag{Ws}{$\OmegaS$}
\psfrag{Wl(t)}{$\OmegaFt$}
\psfrag{Ws(t)}{$\OmegaSt$}
\psfrag{G}{$\varGamma$}
\psfrag{Gsl}{$\varGamma_\text{\!\sc fs}$}
\psfrag{G(t)}{$\varGamma_\text{\!\sc fs}(t)$}
\psfrag{xi}{$\bm\xi(t)$}
\psfrag{x}{$\xx$}
\psfrag{X}{$\XX$}
\psfrag{t6}{$m_\text{\sc s}(0)$}
\hspace*{-4em}\includegraphics[width=37em]{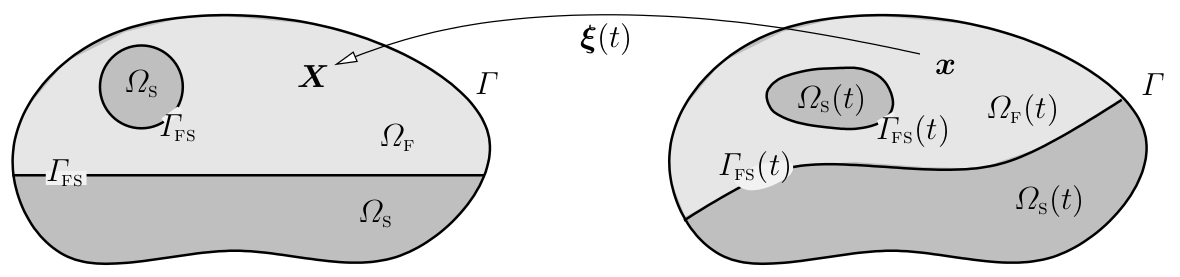}
\end{my-picture}
\nopagebreak
\\
{\small\sl\hspace*{1em}Fig.~\ref{fig2}:~\begin{minipage}[t]{33em}
A reference configuration (left) and a deformed configuration at time $t$
(right) with the shape of $\varOmega$ fixed but the displacement
in the tangential direction along $\varGamma$ possible.
\end{minipage}
}
\end{center}
Specifically, we prescribe $1+(d{-}1)+d=2d$ boundary conditions:
\begin{subequations}\label{BC-Navier}
\begin{align}\label{BC-1-Navier}
&\vv{\cdot}\nn=0,\ \ \ \ 
\\&\label{BC-2-Navier}\!\!\big[\big(\TTRxi(\FF)\,{+}\,\DISRxi(\ee(\vv))
     \,{-}\,{\rm div}\,\mathcal{H}(\nabla\ee(\vv))\big)\nn
{-}\divG\big(\mathcal{H}(\nabla\ee(\vv))\nn\big)\big]_\text{\sc t}^{}\!
+\nu_\flat^{}\vv=\TRACTION,
\\&\label{BC-3-Navier}
\Nabla\ee(\vv){:}(\nn{\otimes}\nn)={\bm0}\,,
\end{align}\end{subequations}
with $\nu_\flat^{}>0$ a boundary viscosity coefficient
while $[\,\cdot\,]_\text{\sc t}^{}$ denotes the tangential part of a vector.
The condition \eq{BC-2-Navier} involves a boundary ``friction'' $\nu_\flat^{}$
and, together with \eq{BC-1-Navier}, forms the {\it Navier boundary conditions}
largely used in fluid dynamics. The mechanical energy dissipation balance
\eq{mech-engr} extends to
\begin{align}\nonumber
  \hspace*{0em}\frac{\d}{\d t}
\int_\varOmega\!\!\!\!
\linesunder{\frac{\phiRxi(\FF)}{\det\FF}}{stored}{energy}\!\!\!\d\xx
+\!\int_\varOmega\!\!\!\!\lineunder{\DISRxi(\ee(\vv)){:}\ee(\vv)
+\nu|\Nabla\ee(\vv)|^s_{_{_{}}}\!}{bulk dissipation rate}\!\!\d\xx
+\!\int_\varGamma\!\!\!\!\!\!\!\!\!\!\linesunder{\nu_\flat^{}|\vv|^2_{_{_{}}}\!}{boundary}{dissipation rate}\!\!\!\!\!\!\!\!\!\!\d S\
\\[-.4em]\nonumber
=\int_\varOmega\!\!\!\!\!\!\!\linesunder{\frac{\rhoRxi\,\GRAVITY{\cdot}\vv}{\det\FF}}{power of}{gravity field}\!\!\!\!\!\!\d\xx
+\int_\varGamma\!\!\!\!\!\!\linesunder{\!\!\TRACTION{\cdot}\vv_{_{_{_{}}}}\!\!}{power of}{traction}\!\!\!\!\!\d S\,.
\end{align}
\end{remark}

\begin{remark}
[{\sl Making the boundary $\varGamma$ free: a  fictitious-domain or  sticky-air approach}]
\label{rem-BC++}\upshape
In some applications, an entirely free $\varGamma$ is urgently desirable.
This seems possible only by an ``engineering'' approach, embedding $\varOmega$
into a fictitious bigger domain containing a very soft material (fluid)
 with a very low viscosity 
and having with a fixed boundary to comply the usual non-penetration
requirement of the outer boundary $\varGamma$  and also to the analysis
presented above, see e.g.\ \cite{PatRae16FEVS}. In geophysical
modelling, such rather rough (although numerical simple and efficient)
trick is called a ``sticky-air'' approach, cf.\ e.g.\ \cite{CSGD12CNST}.
\end{remark}

\begin{remark}[{\sl Transport by non-Lipschitz velocity field}]\label{rem-reglar-transport}\upshape
The multipolar viscosity ensures the velocity field regular, namely
$\vv\in L^s(I;W^{2,s}(\varOmega;\R^d))$, so that, in particular,
it is valued in $W^{1,\infty}(\varOmega;\R^d)$ when $r>d$ is assumed, i.e.\ the
velocity field is Lipschitz continuous in space integrably in time. Since the
seminal paper by R.\,DiPerna and P.L.\,Lions \cite{DiPLio89ODET}, it is well
understood that transport by non-Lipschitz velocity fields can (and even, in
general, must \cite{AlCrMa19LRCE}) lead to development of singularities and
is very difficult, as documented in hundreds of subsequent articles, as in
particular \cite{Ambr08TECP,BahChe94ETRC,CheLer95FCVN,CriDeL08ERRD,Desj97LTEI}. 
Relaxing such higher-gradient multipolar viscosity is therefore great but surely
very difficult challenge. 
\end{remark}

\bigskip

{\small

\baselineskip=12pt

\noindent{\it Acknowledgments.}
Inspiring discussions with Sebastian Schwarzacher during
author's stay at University of Vienna in 2020 are deeply acknowledged.
Also the support from M\v SMT \v CR (Ministry of Education of the
Czech Republic) project CZ.02.1.01/0.0/0.0/15-003/0000493,
and from the institutional support RVO:61388998 (\v CR)  is 
acknowledged.

} 

\bigskip

\noindent
Mathematical Institute, Faculty of Math. \& Phys., Charles University,\\
Sokolovsk\'a 83, CZ-186~75~Praha~8,  Czech Republic,\\[.3em]
and\\[.3em]
Institute of Thermomechanics, Czech Academy of Sciences,\\
Dolej\v skova 5, CZ-18200~Praha~8, Czech Republic\\
email: ${\texttt{tomas.roubicek@mff.cuni.cz}}$

\end{document}